\documentclass[12pt]{article}
\usepackage{amsmath,amssymb,latexsym}
\usepackage{graphicx}
\usepackage{epstopdf}

\def\pf{{\bf Proof }}

\begin{document}
\title{Ruled nodal surfaces of Laplace eigenfunctions and injectivity sets for the
spherical mean Radon transform in $\mathbb R^3.$}
\author{Mark L.~Agranovsky}
\maketitle
\begin{center}
Bar-Ilan University
\end{center}
\newcommand{\bt}{\begin{Theorem}}
\newcommand{\et}{\end{Theorem}}
\newcommand{\bi}{\begin{itemize}}
\newcommand{\ei}{\end{itemize}}
\newcommand{\bea}{\begin{eqnarray}}
\newcommand{\eea}{\end{eqnarray}}
\newtheorem{Definition}{Definition}[section]
\newtheorem{Theorem}[Definition]{Theorem}
\newtheorem{Lemma}[Definition]{Lemma}
\newtheorem{Exercise}{\sc Exercise}[section]
\newtheorem{Proposition}[Definition]{Proposition}
\newtheorem{Corollary}[Definition]{Corollary}
\newtheorem{Problem}[Definition]{Problem}
\newtheorem{Example}[Definition]{Example}
\newtheorem{Remark}[Definition]{Remark}
\newtheorem{Remarks}[Definition]{Remarks}
\newtheorem{Question}[Definition]{QuesMotition}
\newtheorem{Statement}[Definition]{Statement}
\newtheorem{Conjecture}[Definition]{Conjecture}
\newcommand{\be}{\begin{equation}}
\newcommand{\ee}{\end{equation}}
\def\pf{{\bf Proof }}
\vskip.25in

\setcounter{equation}{0}
 \begin{abstract}
It is proved that if a Paley-Wiener family of eigenfunctions of the Laplace operator in
$\mathbb R^3$ vanishes on a real-analytically ruled two-dimensional surface $S \subset
\mathbb R^3$ then $S$ is a union of cones, each of which is contained in a translate of
the zero set of a nonzero harmonic homogeneous polynomial. If $S$ is an immersed $C^1-$
manifold then $S$ is a Coxeter system of planes. Full description of common nodal sets of
Laplace spectra of convexly supported distributions is given. In equivalent terms, the
result describes ruled injectivity sets for the spherical mean transform and confirms,
for the case of ruled surfaces in $\mathbb R^3,$ a conjecture from \cite{AQ1}.
\end{abstract}



\section{Introduction}
Nodal sets are zeros of the Laplace eigenfunctions. They play an important role in
understanding of the wave propagation.

The geometry of a single nodal set can be very complicated  and hardly can be well
understood. On the other hand, simultaneous vanishing of large families of eigenfunctions
on large sets occurs rarely and hence it is naturally to expect that common nodal sets
in that case should be pretty special and have a simple geometry.

Bourgain and Rudnick \cite{BR} obtained a result of such type for two-dimensional torus
$T^2.$ They proved that only geodesics can serve common nodal curves for infinitely many
Laplace eigenfunctions on $T^2.$ For tori in high dimensions, they proved that
Gauss-Kronecker curvature of the common nodal hypersurfaces must be zero. Analogous
question for the sphere in the Euclidean space is still open.

In this article, we address to the similar questions for Euclidean spaces. The case of
$\mathbb R^2$ was studied in \cite{AQ1}, in equivalent terms of injectivity sets for the
spherical mean Radon transform. Translated back to the language of nodal sets, the result
of \cite{AQ1} says that one-dimensional parts of common nodal sets of large families of
eigenfunctions are Coxeter system of straight lines in the plane. There, by large
families of eigenfunctions, Laplace spectra of compactly supported
functions were understood.

In the course  of that result, it was conjectured in \cite{AQ1} that in higher
dimensions, common nodal surfaces for large families of eigenfunctions (injectivity sets
of the spherical mean transform) are cones - translates of the zero sets of solid
harmonics (harmonic homogeneous polynomials). In this article, we confirm  Conjecture
from \cite{AQ1} for a special case of ruled surfaces in $\mathbb R^3.$ The proof develops ideas from
the article \cite{AQ2} of E.T. Quinto and the author.

Although ruled surfaces (unions of straight lines) are, in a sense, close to cones (union
of straight lines with a common point), proving conical structure of ruled nodal surfaces
in dimensions higher than two was elusive for a long time.

\section{Content}
\begin{itemize}
\item 1. Introduction
\item 2. Content
\item 3. Main results
\item 3.1.  Nodal surfaces version
\item 3.2. Injectivity sets version
\item 4. Background
\item 5. The strategy of the proof of the main result
\item 6. Preliminary observations
\item 7. Local symmetry and antipodal  points
\item 8. Ruled surfaces
\item 8.1. Regularity of the line foliation at smooth points
\item 9. The structure of real analytically ruled algebraic surfaces near singular points
\item 9.1. The outline of the proof
\item 9.2. Preliminary constructions
\item 9.3. Re-scaling
\item 9.4. Re-parametrization
\item 9.5. The case of odd $m$
\item 9.6. The case of even $m$
\item 9.7. End of the proof of Theorem \ref{T:conical}
\item 10.  Irreducible case. Proof of Theorem \ref{T:Main1}
\item 10.1. Extremal ruling lines and antipodal points
\item 10.2. End of the proof of Theorem \ref{T:Main1}
\item 11. Reducible case. Proof of Theorem \ref{T:mainmain1}
\item 12. Coxeter systems of planes. Proof of Theorem \ref{T:mainmain2}
\item 13. Proof of Theorem \ref{T:main_conv_nodal} (the case of convexly supported generating function)
\item 14. Concluding remarks
\item 15. References
\end{itemize}

\section{Main results}

We will formulate the main results of this article in two equivalent terms: 1) on the
language of nodal surfaces and 2) on the language of of injectivity sets.

We start with the nodal surfaces version.
\subsection{Nodal surfaces version}
Let $\varphi_{\lambda}, \ \lambda >0, $ be a family of eigenfunctions of the Laplace
operator $\Delta$  in $\mathbb R^3.$ More precisely, each function $\varphi_{\lambda}$ is
a solution of the Helmholtz equation
$$\Delta \varphi_{\lambda}=-\lambda^2 \varphi_{\lambda}.$$
In particular, $\varphi_{\lambda}$ can be identically zero function.
\begin{Definition} The family $\varphi_{\lambda}$ is a {\bf Paley-Wiener} family
if it can be extended in the complex plane $\lambda \in \mathbb C$ as an even entire
function, satisfying the growth condition
$$|\varphi_{\lambda}(x)| \leq C (1+|\lambda|)^N e^{(R+|x|) |Im \lambda|}.$$
for some positive constants $C,R$ and for some natural $N.$
\end{Definition}

By cone in $\mathbb R^d$, we understand union of straight lines having a common point-the
vertex of the cone. We call a cone $C$  {\bf harmonic cone} if there exists a nonzero
harmonic homogeneous polynomial (solid harmonic) $h$  and a vector $a$ such that
$$C \subset a+h^{-1}(0).$$

\begin{Definition} \label{D:ruled_surface} Let $S$ be a surface in $\mathbb R^3.$ We call $S$ {\bf
irreducible real analytically ruled surface} if
\begin{enumerate}
\item There exists a  closed continuous curve $\gamma \subset \mathbb R^3$ such that
$S$ is the union of straight lines, $ S=\cup_{a \in \gamma} L_a,$ passing through points
$a\in \gamma.$
\item Locally, $S$ is the image of the (parametrizing) mapping
$$(-1,1) \times \mathbb R \ni  (t,\lambda) \to u(t,\lambda) =u(t)+\lambda e(t),$$ where
$I \ni t \to (u(t), e(t)) \in \mathbb R^3 \times \mathbb R^3$ are real analytic maps and $|e(t)|=1.$
\end{enumerate}
The curve $\gamma$ is called the {\bf base curve}, the vector
$e(t)$-directional vector, the straight lines $L_t=\{u(t)+\lambda e(t), \ \lambda \in
\mathbb R \}$ are called rulings, or ruling or generating lines.
{\bf Real analytically
ruled surface} are, by definition, finite unions of irreducible those.
\end{Definition}

\begin{Remark}
\begin{enumerate}
\item
The parametrizing mapping $u(t,\lambda)$ is nor necessarily defines a parametrization of $S$ as a manifold,
since the regularity condition is not required.
\item
Real analytically ruled surfaces are not necessarily  everywhere real analytic, and
even differentiable. For example, a cone which is not a plane is not
differentiable at its vertex.
\end{enumerate}
\end{Remark}

Now we are ready to formulate the main results of this article.

\begin{Theorem} \label{T:Main1}
Let $S$ be an irreducible  real-analytically ruled surface with no parallel generating
lines, then $S$ is the common nodal set for a Paley-Wiener family if and only if $S$ is a
harmonic cone.
\end{Theorem}

In the reducible case, we have
\begin{Theorem} \label{T:mainmain1} Let $S$ be a real-analytically ruled surface in
$\mathbb R^3,$ with no parallel generating lines. If $S$ is the common nodal for a
Paley-Wiener family of eigenfunctions then $S$ is the union of finite number  of harmonic
cones, $S=\cup_{j=1}^N C_j$ such that for any $1 \leq i < j \leq N$ the intersection $C_i
\cap C_j \neq \emptyset$ and  one the two cases are possible:
\begin{enumerate}
\item $C_i \cap C_j$ is the vertex of one of the cones $C_i,C_j,$
\item $C_i \cap C_j$ is transversal and is an unbounded curve.
\end{enumerate}
\end{Theorem}

Conjecture from \cite{AQ1} (see section \ref{S:background} for the details) claims that,
in fact, the vertices of the cones $C_i$ coincide and therefore $S$ is a single cone.
However, we are not able to prove that at the moment.

\begin{Definition} The union $\Sigma=\cup_{j=1}^N \Pi_j$ of $N$ hyperplanes
in $\mathbb R^d$ having a common point is called {\bf Coxeter system} if $\Sigma$ is
invariant with respect to  all the reflections around the planes $\Pi_j, \ j=1,...,N.$
\end{Definition}

Notice that Coxeter systems are harmonic cones, i.e., are, up to translations, zero sets
of solid harmonics.
\begin{Theorem} \label{T:mainmain2} If in Theorem \ref{T:mainmain1} $S$ is an immersed
$C^1-$ surface then $S$ is a Coxeter system.
\end{Theorem}

Remind that immersed $C^1-$ surface is the image of a two-dimensional $C^1-$ manifold
under a $C^1-$  mapping with non-degenerating differential.

Finally, we will formulate  one more result about common nodal surfaces for special
Paley-Wiener families of eigenfunctions: spectral projections of  {\it convexly
supported} distributions:
\begin{Theorem}\label{T:main_conv_nodal}  Let $f \in D^{\prime}_{comp}(\mathbb R^3)$
be a nonzero compactly supported distribution or continuous function  and
$$f=\int_{0}^{\infty}\varphi_{\lambda}d\lambda$$ be the Laplace spectral decomposition of
$f$ (\cite{Strichartz}). Assume that the boundary of the unbounded connected component of
$\mathbb R^3 \setminus supp f$ is a real analytic strictly convex closed surface. If
$$N=\cap_{\lambda>0}\varphi_{\lambda}^{-1}(0)$$ then $N=S \cup V$ where either $V= \emptyset$ or
$V$ is an algebraic variety of $dim V \leq 1$ and either $S = \emptyset$ or $S$  is one of the  three
surfaces:
\begin{enumerate}
\item $S$ is a harmonic cone.
\item $S$ is the union of two harmonic cones, $S=C_1 \cup C_2$ such that either $C_1 \cap C_2 =\{b_1\}$ or $C_1 \cap
C_2=\{b_2\}.$ where $b_1,b_2$ are the vertices of the corresponding cones.
\item $S$ is the union of three harmonic cones, $S=C_1 \cup C_2 \cup C_3,$ with the
vertices $b_1, b_2, b_3,$ correspondingly, such that either
$$C_1 \cap C_2=\{b_1\}, \ C_2 \cap C_3=\{b_2\}, \ C_3 \cap C_1 =\{b_3\}$$ or
$$C_1 \cap C_2=\{b_2\}, \ C_2 \cap C_3=\{b_3\}, \ C_3 \cap C_1 =\{b_1\}.$$
\end{enumerate}
\end{Theorem}

We conjecture that, in fact,   $b_1=b_2=b_3$ which would lead to confirming Conjecture \ref{C:Conjecture} in a
more complete form.

\subsection{Injectivity sets version}\label{S:IS}

The spherical mean Radon transform  is defined  as the mean value
$$Rf(x,t)=\int\limits_{|\theta|=1}f(x+t\theta)dA(\theta)$$
of $f$ over the sphere $S(x,t)$ centered at $x \in \mathbb R^d$ of radius $t >0.$ Here
$dA$ is the normalized area measure on the unit sphere $\{|\theta|=1\}$ in $\mathbb R^d.$

The operator $R$ can be extended to distributions $f \in D^{\prime}(\mathbb R^d).$
Namely, for each vector $a \in \mathbb R^d$ define the averaging operator
$$R_a\psi(x):=\int_{SO(n)} \psi(a+\omega(x-a))d\omega,$$
where $d\omega$ is the normalized Haar measure on the orthogonal group $SO(n).$ The
relation between this averaging operator and the operator $R$ is given by
$$(R_a\psi)(x)=R\psi(a,|x-a|).$$
Now, if $f \in D^{\prime}(\mathbb R^d)$ and $a \in \mathbb R^d,$ then we define the new
distribution $R_af$ by the following action on test-functions $\psi$:
\begin{equation}\label{E:R_a}
 \langle R_a f,\psi \rangle = \langle f, R_a\psi \rangle.
\end{equation}
It is easy to see that  this definition is consistent with the definition of the action of the operator
$R_a$ on functions.

Denote $R_S$ the restriction of the transform $R$ on the set $S \times (0,\infty):$
$$ R_S: C_{comp}(\mathbb R^d) \ni f \to Rf\vert_{S \times (0,\infty)}.$$

\begin{Definition}\label{D:inj_sets}  We call a
set $S \subset \mathbb R^d$ {\bf injectivity set} if given a distribution $f \in
D^{\prime}_{comp}(\mathbb R^d)$ such that $R_a f=0$ for all $a \in S$ then $f=0.$
Equivalently, $S$ is injectivity set if the operator $R_S$ is injective,
 i.e. for every function $f \in C_{comp}(\mathbb R^d)$
$$Rf(x,t)=0 \ \mbox{for all} \ \ x \in S \ \  \mbox{implies} \ \  f=0.$$
\end{Definition}

The equivalence of definition for functions and distributions can be easily proved by
convolving distributions with radial smooth functions.

Spherical mean Radon transform \footnote{We refer to Radon transform because the operator
$R$ is defined on complexes of spheres with restricted centers and of  arbitrary radii.
Such varieties are analogous to  varieties of planes with restricted set of normal
vectors and arbitrary distances to the origin which are natural in the study of the plane
Radon transform.} plays an important role in applications, namely, in thermo- and
photoacoustic tomography (cf. \cite{KUCHMENT}), which is used in the medical imaging \cite{KUCH}.
The mathematical problem behind `that is to recover $f$ from the data
$Rf(x,t), \ x \in S, \ t >0.$
The uniqueness of the recovery is equivalent to the injectivity of the operator $R_S$ and
therefore the first question to be answered is to understand for what {\bf observation
surfaces} $S$ the operator $R_S$ is injective, i.e., to understand the injectivity sets.
Of course, the case $d=3$ is most important from the point of view of the applications.

\begin{Definition} Let $\{\varphi_{\lambda}\}_{\lambda >0,}$ be a measurable family of
Laplace eigenfunction: $(\Delta +\lambda^2)\varphi_{\lambda}=0$ in $\mathbb R^d.$ We will
call the function
\begin{equation}\label{E:f}
f(x)=\int\limits_0^{\infty} \varphi_{\lambda}(x)\lambda^{d-1} d\lambda
\end{equation}
{\bf generating function}, assuming that the integral converges ( which can be achieved
by a proper normalization $\varphi_{\lambda} \to c(\lambda) \varphi_{\lambda}.$) The
family $\varphi_{\lambda}$ is called a {\bf Laplace spectral decomposition} of $f.$
\end{Definition}
The definition can be extended to distributions $f \in D^{\prime}(\mathbb R^d)$ if to
understand the spectral decomposition of $f$ in the distributional sense.

The weight factor $\lambda^{d-1}$ is added for convenience and can be omitted by
including it to $\varphi_{\lambda}.$

The link between common nodal sets  and injectivity sets in the question is very simple: they just coincide (see Proposition \ref{P:observations}).

Let us briefly explain this relation. It is proved in (\cite{Strichartz}, Theorem 3.10) that

{\it the family $\varphi_{\lambda}$ of
eigenfunctions in $\mathbb R^d$ is Paley-Wiener if (and if and only if , when $d$ is odd)
the integral (\ref{E:f}) defines a compactly supported distribution $f \in D^{\prime}(\mathbb R^d).$}

The spectral decomposition $\{\varphi_{\lambda}\}$ can be recovered from  the generating distribution $f$ by means of
the convolutions
\begin{equation}\label{E:spect}
\lambda^{d-1} \varphi_{\lambda}=\j^{\lambda}_{\frac{d-2}{2}} * f
\end{equation}
of $f$ with the normalized Bessel function
$$j^{\lambda}_{\frac{d-2}{2}}(x)=(2\pi)^{-\frac{d}{2}} \frac{J_{\frac{d-2}{2}}(|\lambda
x|)}{(|\lambda x|)^{\frac{d-2}{2}}}.$$
It follows that
$S \subset \cap_{\lambda>0} \varphi_{\lambda}^{-1}(0)=0$ if and only if
$Rf\vert_{S \times (0,\infty)}=0.$

Remind that the condition
$Rf\vert_{S \times (0,\infty)}=0$ for $f \in D^{\prime}(\mathbb R^d)$
means that the average distribution $R_a f,$ defined in (\ref{E:R_a}), is the zero
distribution: $R_a f=0$ for all $a \in S.$

Thus, we have
\begin{Proposition}
A set $S \subset \mathbb R^d$ serves a common nodal set for a nontrivial family
$\{\varphi_{\lambda}\}$ if and only if $Rf\vert_{S \times (0,\infty)}=0$ for some nonzero
compactly supported distribution (or continuous function) $f$, i.e., if and only if $S$
fails to be a set of injectivity for the spherical mean Radon transform $R.$
\end{Proposition}
Using that equivalence, we can reformulate  Theorems \ref{T:Main1} and  \ref{T:mainmain1}  in
the equivalent form:

\begin{Theorem}\label{T:reform1} Let $S$ be a real-analytically ruled surface in $\mathbb R^3.$
If $S$ fails to be an injectivity set then $S$ is one of the surfaces enlisted in Theorem
\ref{T:mainmain1}. If $S$ is irreducible (see Definition \ref{D:ruled_surface}) then $S$
fails to be an injectivity set if and only if $S$ is a harmonic cone.
\end{Theorem}

The following theorem is a translation, on the injectivity sets language, of Theorem
\ref{T:main_conv_nodal}. It follows from Theorem \ref{T:mainmain1} and
\cite{AQ1}, \cite{AKuch} and is an equivalent version of Theorem \ref{T:main_conv_nodal}.
Here the condition refers to the geometric shape of the support of the generating
distribution.

\begin{Theorem}\label{T:main_conv_inj}  Let $f \in D^{\prime}_{comp}(\mathbb R^3)$
be nonzero compactly supported distribution or continuous function. Assume that the
boundary of the unbounded connected component of $\mathbb R^3 \setminus supp f$ is a
real analytic strictly convex closed surface. If $Rf(x,t)=0$ for all $x \in S$ and $t>0$
then $S$ is one of the surfaces enlisted in Theorem \ref{T:mainmain1}.
\end{Theorem}
The proof of Theorems \ref{T:main_conv_nodal} and \ref{T:main_conv_inj} is based on
Theorem \ref{T:mainmain1} and the results of \cite{AQ1},\cite{AKuch} (Theorem
\ref{T:ruled}, the next section) about ruled structure of observation surfaces for
convexly supported functions.

\section{Background}\label{S:background}

In dimension $d=2$, the problem of describing injectivity sets was completely solved in
[AQ1]. Let us formulate the result. Denote
$$\Sigma_N=(tcos\ k\frac{\pi}{N},t sin \ k\frac{\pi}{N}), \ k=0,1,...,N-1, \  -\infty < t <\infty, $$
the (Coxeter) system of $N$ straight lines passing through the origin and having equal angles
between the adjacent lines.

\begin{Theorem} \label{T:AQ1} [\cite{AQ1}] A set $S \subset \mathbb R^2$ is a set of injectivity if
and only if $S$ is contained in no set of the form $(a+ \omega (\Sigma_N)) \cup V ,$
where $a \in \mathbb R^2, \ \ \omega$ is a rotation in the plane and $V$ is a finite set,
invariant under reflections around the lines from the Coxeter system $a+ \omega
(\Sigma_N)$.
\end{Theorem}

Observe that the Coxeter system  $\omega(\Sigma_N)$ coincides with the zero set of the
polynomial $h(x,y)=Im (e^{i\varphi}(x+iy)^N),$ where $\omega$ is the rotation for the angle
$\varphi.$ The polynomial $h(x,y)$ represents the general form of harmonic homogeneous polynomial in the plane.
That observation  gives rise to the following conjecture about how injectivity sets look
like in arbitrary dimension.

\begin{Conjecture} \label{C:Conjecture} [\cite{AQ1}] Let $S \subset \mathbb R^d$ fails to be an
injectivity set Then $S \subset (a+h^{-1}(0)) \cup V,$ where $h$ is a harmonic
homogeneous polynomial (spatial harmonic) and  $V$ is an algebraic variety in $\mathbb
R^d$ of dimension $dim V
 \leq d-2.$
\end{Conjecture}

Since in odd dimensions, as it was mentioned in subsection \ref{S:IS}, non-injectivity
sets are precisely common nodal sets of Paley-Wiener families, Conjecture \ref{C:Conjecture} can be reformulated as following:

\begin{Conjecture} \label{C:conj2} A set $S \subset \mathbb R^d, \ d$ is odd, is a
common nodal set for a  Paley-Wiener family of Laplace eigenfunctions if and only  $S
\subset (a+h^{-1}(0)) \cup V,$ where the vector $a,$ the variety $V$ and the polynomial
$h$ are as in Conjecture \ref{C:Conjecture}.
\end{Conjecture}

\begin{Remark}
A partial case of non-injectivity sets in Conjecture are Coxeter systems of hyperplanes.
They are arrangements of $N$ hyperplanes with a common point, invariant under reflections
around each the hyperplane from the system. The Coxeter systems correspond to the case of
completely reducible harmonic homogeneous polynomials $h,$ i.e., those represented as
products $$ h=l_1 \cdots l_N$$ of $N=deg h$ linear forms.
 \end{Remark}

Here are some evidences for Conjecture  \ref{C:Conjecture} (see \cite{AParab}):

\begin{itemize}
\item  Any harmonic cone is a non-injectivity set, i.e.,  if  $h$ is a non-zero harmonic
homogeneous polynomial, then $S:=h^{-1}(0)$ is a
non-injectivity set. Namely, define $f(x):=\alpha(r)h(\theta)$ where $r,\theta$ are the
spherical coordinates:$x=r\theta, |\theta|=1$ and $\alpha (r)$ is a non-zero smooth even
compactly supported function on $\mathbb R.$ It is an easy exercise to prove that
$Rf(x,t)=0$ for all $x \in S, \ t>0.$
\item If $V$ is an algebraic variety of $dim V \leq d-2$ then there exists a nonzero
$f \in C_{comp}(\mathbb R^d)$ such that $Rf (x,t)=0$ for all $(x,t) \in V \times
(0,\infty).$ (\cite{AParab}, Theorem 3.2).
\end{itemize}

So far, only partial results towards Conjecture \ref{C:Conjecture} are obtained
\cite{AQ2},  \cite{AKuch}, \cite{AVZ}.

It was proved in  \cite{AVZ} that among cones only zero sets of spatial
harmonics fail to be injectivity sets. Therefore, the main difficulty in proving
Conjecture \ref{C:Conjecture} is checking that non-injectivity sets are necessarily {\it
cones}.

The following two results can be considered as  certain steps in that direction:

\begin{Theorem} [\cite{AQ3}]
Let $f$ be a compactly supported continuous function or distribution in $\mathbb R^d.$
Assume that $supp f$ is the union of disjoint balls or $supp f$ is finite. If $S \subset
\mathbb R^d$ and $R_Sf=0$ then $S \subset (a+h^{-1}(0)) \cup V,$ where $a \in \mathbb
R^d, \ h$ is a nonzero harmonic homogeneous polynomial and $V$ is an algebraic variety of
$dim V \leq d-2.$
\end{Theorem}

The next result deals with functions with convex compact supports and can be viewed as a motivation for Theorems
\ref{T:main_conv_nodal} and \ref{T:main_conv_inj}.
\begin{Theorem} \label{T:ruled}[ \cite{AKuch}, \cite{AQ2}].
Let $f \in C_{comp}(\mathbb R^d)$ be a compactly supported function.
Suppose that the outer boundary $\Gamma=\partial {supp f}$ is convex. If  $S \subset
\mathbb R^d$ is such that $Rf\vert_{S \times (0,\infty)}=0$ then $S$ is ruled, i.e., $S$
is the union of straight lines. Moreover, the ruling lines intersect $\Gamma$
orthogonally at each point where $S$ is differentiable.
\end{Theorem}

By outer boundary $\partial supp f$ we understand the boundary $\partial (\mathbb R^d \setminus supp f)_{\infty}$
of the unbounded connected component of the complement.
\begin{Remark}  In fact,  the ruled structure of $S$ was established in  \cite{AKuch}
under much milder conditions for $\Gamma$ for example, under assumption of $C^2$
smoothness of $\Gamma.$ However,  in the proofs of Theorems \ref{T:main_conv_nodal} and
\ref{T:main_conv_inj},  we will use the weaker version, Theorem \ref{T:ruled}, because
some additional properties delivered by the convexity of support will be exploited.
\end{Remark}

\section {The strategy of the proof of the main result}

The main result of this article is Theorem \ref{T:Main1}. Theorem \ref{T:mainmain1} is
deduced from Theorem \ref{T:Main1}, Theorems \ref{T:mainmain2}, \ref{T:main_conv_nodal}
follow from Theorems \ref{T:Main1} and \ref{T:mainmain1}. All the theorems can be viewed
as results towards  proving Conjecture \ref{C:Conjecture}-\ref{C:conj2}.

The proof of Theorem \ref{T:Main1} falls apart into several steps:

{\bf Step 1.} First, we prove that the common nodal surface $S$ for Paley-Wiener family
is algebraic and lies in the zero set of a nontrivial harmonic polynomial. In a different
setting, that fact was first observed in \cite{LP} (see also \cite{AQ2}).

{\bf Step 2.} Next, we formulate local symmetry property, which is based on the results of
\cite{AQ1},\cite{SU} about cancelation of analytic wave front sets.  The corollary of
that property says is that any surface $S$ having a pair of antipodal points-points of
smoothness, such that the segment joining them is orthogonal to the surface, fails to be
a common nodal surface for a Paley-Wiener family.

{\bf Step 3.} Assuming that $S$ is not a cone and using compactness argument we find two
generating (ruling) straight lines on $S$ with the maximal distance between them. Then we
pick two closest points $a, b \in S$ on those extremal lines. If those extremal points
$a,b$ are regular then the previous step implies that $S$ cannot be nodal. Otherwise, one
of the extremal points is singular and we encounter the problem of characterization of
singularities of algebraic real analytically ruled surfaces in $\mathbb R^3.$

{\bf Step 4.}  We obtain the required characterization of the singularities (Theorem \ref{T:conical}), which is a
key ingredient of the proof of the main result.

{\bf Step 5.} The final arguments are as follows.  Theorem \ref{T:conical} claims that singular points
are conical or of cuspidal type. The corollary \ref{C:harm} is that either the
irreducible ruled algebraic surface $S$ is a cone or it is a uniqueness set of harmonic
polynomials. However,the latter option is ruled out (Step 1). Therefore,  we conclude
that $S$ is a cone (in the irreducible case) or a union of cones (in the reducible case).
Finally, the proof that the cones are harmonic easily follows by homogenization of
harmonic polynomial vanishing on $S$ (obtained on Step 1). This completes the proof.

\begin{Remark} Essentially, steps 1-3 were presented in the  \cite{AQ2}.
It was proved there that if the extremal points (Step 3) are regular then the surface is
an injectivity set (not nodal). The description of singular points obtained in Theorem
\ref{T:conical} allowed us to further develop the idea of \cite{AQ2} and push forward
proving the conical structure of the nodal ruled surfaces, which is the main
result of this article.
\end{Remark}

\section{Preliminary observations}

In this section, we  briefly present auxiliary facts that we will need in the sequel. Most of
them are exposed  in \cite{AQ1}. It will be convenient to combine those facts in one
proposition:
\begin{Proposition}\label{P:observations} Let $\Phi=\{\varphi_{\lambda}, \ \lambda >0,\}$
be a nonzero family of Laplace  eigenfunctions in $\mathbb R^d$ with compactly supported
generating distribution $f \in D^{\prime}(\mathbb R^d)$ i.e.,
$$f=\int_{0}^{\infty}\varphi_{\lambda}d\lambda.$$
We omit the factor $\lambda^{d-1}$ by including it into $\varphi_{\lambda}.$ Clearly, this does not effect on the zero sets of $\varphi_{\lambda}.$

 Denote
$$N_f =\{x \in \mathbb R^d: Rf(x,t)=0, \ \forall   (x,t) \in S \times (0,\infty) \}$$
and
$$N(\Phi)=\cap_{\lambda >0}\varphi_{\lambda}^{-1}(0).$$
Then
\begin{enumerate}
\item
$N_f=N(\Phi)$ and therefore common nodal sets and non-injectivity sets are the same.
\item
The set $N(\Phi)$ is algebraic and has the form
$$N(\Phi)=S \cup V,$$
where $S=\emptyset$ or $S$ is a real algebraic hypersurface: $S=Q^{-1}(0),$ where $Q$ is
a nonzero real polynomial, and $V$ is an algebraic variety of $dim V \leq d-2$ (maybe,
empty as well).
\item
There is a nonzero real harmonic polynomial $H$ vanishing on $S,$ i.e. $S \subset
H^{-1}(0).$
\end{enumerate}
\end{Proposition}

{\bf Proof}
\begin{enumerate}
\item
 We have $f=\int_0^{\infty} \varphi_{\lambda} d\lambda,$  where the equality is
understood in the distributional sense.

Then  for any test-function $\psi$ and for any $a \in \mathbb R^d:$
$$ \langle R_af,\psi \rangle = \langle f, R_a\psi \rangle =\int_{\mathbb R} \langle \varphi_{\lambda}(x), R_a\psi \rangle d\lambda.$$

Further,
$$\langle \varphi_{\lambda}, R_a\psi \rangle =\int_{\mathbb R^d}\int_{SO(n)}\varphi_{\lambda}(x)
\psi(a+\omega(x-a))dxd\omega.$$
Change of variables
$y=a+\omega(x-a)$ yields
$$\langle \varphi_{\lambda}(x), R_a\psi \rangle = \langle R_a\varphi_{\lambda},\psi \rangle .$$
The Laplace eigenfunctions are also eigenfunctions of the averaging operator $R_a:$
$$(R_a \psi_{\lambda})(x)=j_{\frac{d-2}{2}}(\lambda |a|)\psi_{\lambda}(x).$$
Therefore, we have
$$\langle R_a f,\psi \rangle =\langle \int_0^{\infty} j_{\frac{d-2}{2}}(\lambda |x-a|)\varphi_{\lambda}(a)
d\lambda, \psi(x)\rangle_{x}.$$
We see that $R_af=0$ if and only if
$$\int_0^{\infty} j_{\frac{d-2}{2}}(\lambda |x-a|)\varphi_{\lambda}(a)
d\lambda=0$$ for all $x$. The latter integral equation is satisfied if and only if
$\varphi_{\lambda}(a)=0.$ Thus, $N_f=N(\Phi)$ and the statement 1 is proved.

\item Decompose the (even) normalized  Bessel function $j_{\frac{d-2}{2}}(\lambda t)$ into
power series:

$$j_{\frac{d-2}{2}}(\lambda t)=\sum_{k=0}^{\infty} c_k \lambda^{2k} t^{2k}.$$

Then we have from (\ref{E:spect}):
$$\varphi_{\lambda} (x)=\sum_{k=0}^{\infty} c_k \lambda^{2k} |x|^{2k} * f.$$
Therefore, $x \in N(\Phi)$ if and only if $|x|^{2k} *f=0, \ k=0,1,....$

Notice that
$$Q_k(x)=c_k|x|^{2k} * f=c_k <|x-y|^{2k}, f>,$$
where the right hand side stands for the action of the distribution $f$ with respect to
$y.$ It follows that $Q_k$  is a polynomial, of $deg \ Q_k \leq 2k.$

From (\ref{E:spect})  $\varphi_{\lambda}(x)=0$ is equivalent to $Q_k(0)=0, k=0,1,...$ and
hence
$$N(\Phi)=\cap_{k=0}^{\infty}Q_k^{-1}(0).$$

Denote $Q$ the greatest common divisor (over $\mathbb C$) of $Q_k.$ Then
$$N(\Phi)=(Q^{-1}(0) \cap \mathbb R^d) \cup V,$$
where $V$ is the intersection of $\mathbb R^d$ with the zero varieties of coprime
polynomials and hence $dim_{\mathbb R} V < d-1.$

To complete the proof of the statement 2, we have to show that the polynomial $Q$ has real coefficients.
We will do that at the end of the proof.

\item
 Substitution (\ref{E:spect}) into Helmholtz equation :
$$\Delta \sum_{k=0}^{\infty}  \lambda^{2k} Q_k=
-\lambda ^2 \sum_{k=0}^{\infty}  \lambda^{2k} Q_k$$ yields
$$ \Delta Q_k=-Q_{k-1}.$$

Not all polynomials $Q_k$ are identically zero. Indeed, suppose that $Q_k=c_k |x|^{2k} *f
\equiv 0$ for all $k=0,1,...$ Since $f$ has compact support and the linear combinations
of the polynomials $|y|^{2k}$ approximate, in the $C^{\infty}$ topology on compact sets,
any radial smooth function $\alpha (|y|^2)$, we have $\alpha *f \equiv 0.$ Taking Fourier
transform , we obtain $\hat \alpha  \hat f \equiv 0$ which implies $\hat f=0$ due to the
arbitrariness of the radial function $\alpha.$  Then $f=0$ which is not true.

 Let  $k=k_0$ be the minimal $k$ such that $Q_k \neq 0$ and denote
$$H=Q_{k_0}.$$
Then
$$\Delta H=-Q_{k_0-1}=0$$
and hence $H$ is harmonic. This proves the statement 3.

It remains to prove that, in fact, $Q$ is a real polynomial, i.e. has the real
coefficients.  To this end, we first will prove the third statement.

Let
$$H=H_1 \cdots H_q$$
be the decomposition into irreducible , over $\mathbb C,$ polynomials. Let us prove that
all polynomials $H_i$ are real.

Consider the operation of complex conjugation of coefficients:
$$H^*(z)=\overline{H(\overline z)}, \ z \in \mathbb C^d.$$
Since $H=Q_{k_0}$ has real coefficients, we have
$$H^*=H_1^* \cdots H_q^*=H_1 \cdots H_q.$$

Therefore, each $H^*_i$ coincides with some $H_j.$  If for some $i \neq j$ holds
$H_i^*=H_j$ then $H$ is divisible by $H_i H_i^*$ and represents as
$$H=H_i H_i^* R,$$
for some polynomial $R.$ Since in the real space $\mathbb R^d$ we have $H_i^*=\overline
{H_i},$ we have in $\mathbb R^d:$
$$H=|H_i|^2 R.$$
However, Brelot-Choquet theorem \cite{BC} says that no non-negative real polynomial can
divide a real nonzero harmonic polynomial. Therefore, the only possibility is that
$H_i=H_i^*$ for all $i.$  That means that $H_i$ are real polynomials.

The greatest common divisor $Q$ divides $H$ and therefore is a product of some $H_i.$
Since every polynomial  $H_i$ has real coefficients, $Q$ does so.

If $Q$ is constant, i.e., all $Q_k$ are coprime, then $ S=Q^{-1}(0)=\emptyset.$
Otherwise, $S$ is a hypersurface in $\mathbb R^d.$  Indeed, if $dim S <n-1$  then $\mathbb R^d
\setminus Q^{-1}(0)$ is connected and hence everywhere $Q \geq 0$ everywhere  or $Q \leq
0.$  However, this impossible, since Brelot-Choquet theorem states that preserving sign
polynomials cannot divide harmonic polynomials. This completes the proof of Proposition.
\end{enumerate}

\section{Local symmetry and  antipodal points}

\begin{Definition} \label{D:antipodal_points} Let $S \subset \mathbb R^d$ and let
$a,b \in S, \ a \neq b,$ be two distinct points in $S.$ We call $a$ and $b$ {\bf
antipodal points} if
\begin{enumerate}
\item $S$ is a $C^1-$ hypersurface near the points $a,b.$
\item
$ a-b \perp T_a(S), \ a-b \perp T_b(S),$ where $T_a(S), T_b(S)$ are tangent spaces to $S$
at $a$ and $b$ correspondingly.
\end{enumerate}
\end{Definition}
\begin{Theorem}\label{T:antipodal} (\cite{AQ1}, \cite{AQ2}). If $S \subset \mathbb R^d$
has a pair of antipodal points $a,b$ and $S$ is real analytic in neighborhoods of those
points, then $S$ is an injectivity set.
\end{Theorem}

{\bf Example} The hyperboloid $x_1^2+x_2^2-x_3^2=1$ in $\mathbb R^3$ has antipodal points, for example, $(\pm
1,0,0)$ and hence is an injectivity set.

The proof of Theorem \ref{T:antipodal} is based on the following theorem about
certain symmetry of the support of functions with zero spherical means on a surface:

\begin{Theorem} \label{T:local_symm}
(\cite{AQ1}). Let $S$ be a real analytic hypersurface and  $a \in S.$ Let $f
\in C_{comp}(\mathbb R^d)$ be a compactly supported function such that $Rf\vert_{S \times
(0,\infty)}=0.$ Let $x \in supp f$ be a point of local extremum for the distance function
$d(x):=|x-a|$ and denote
$$x^*=x-2 \langle x-a,\nu_a \rangle \nu_a$$
($\nu_a$ is the unit normal vector of $S$ at $a$), the point, symmetric to $x$ with
respect to the tangent plane $T_a(S)$ (mirror point). Then $x^* \in supp f.$
\end{Theorem}

The proof of Theorem \ref{T:local_symm} uses microlocal analysis and results about cancelation of analytic wave front sets at mirror points (\cite{AQ1}, \cite{HZ},\cite{FSU},\cite{SU}).

We are going to exploit Theorem \ref{T:local_symm} for algebraic surfaces $S=Q^{-1}(0),$ where
$Q$ is a real nonconstant polynomial. However, Theorem \ref{T:local_symm} cannot be applied directly as
$S$ is not necessarily everywhere real analytic and, moreover, even differentiable.
Nevertheless, $S$ is real analytic everywhere outside of the critical set
$$crit S:= \{x \in S: \nabla Q(x)=0 \},$$ which is a nowhere dense subset of $S.$
It is enough to establish a local symmetry property, though in a slightly weaker form than in Theorem \ref{T:local_symm}.

Let us introduce some notations and definitions. Given a point $a \in
S$ in a neighborhood of which $S$ is $C^1$ surface we denote
$$\sigma_a: x \to x-2 \langle x-a,\nu_a \rangle \nu_a,$$
the reflection of $\mathbb R^d$ around the tangent plane $T_a(S).$ Here $\nu_a,$ as
above, is the unit normal vector to $S.$

For any $a \in S$ and $r>0$ denote
$$K_{a,r}:=\{x \in supp f:|x-a|=r \}, $$
the intersection of $supp f$ with the sphere $S_r(a)=\{|x-a|=r\}.$

\begin{Theorem}(Local Symmetry Property) \label{T:local_symm_mod}  Let $S \subset \mathbb R^d$ be a  hypersurface, real
analytic except for a nowhere dense subset. Let $f \in C_{comp}(\mathbb R^d)$ be such
that $Rf \vert_{S \times (0,\infty)}=0.$ Let $a \in S$ be a $C^1$ point. Define
$$r=max \{|x-a|: x \in supp f \}.$$ Then $$ \sigma_a(K_{a,r}) \cap supp f \neq \emptyset.$$
\end{Theorem}

\noindent
{\bf Proof} is based on compactness arguments.

Denote for simplicity $K=K_{a,r}, K^*=\sigma_a(K_{a,r}).$ If $E \subset S $ is the set where $S$ is not real analytic , the the point $a$ is a limit point of $S \setminus E$ and hence  we can find a sequence $a_n \in S \setminus E$
such that
$$\lim\limits_{n \to \infty} a_n = a. $$
The surface $S$ is real analytic at any point $a_n$ and the tangent planes
$$T_{a_n}(S) \to T_a(S), \ n \to \infty.$$

Denote
$$r_n=max\{|a_n-x|: x \in supp f \}$$
and let $x_n \in supp f$ be such that
$$|a_n-x_n|=r_n.$$
By the construction, for all $x \in supp f$ holds
$$|a_n-x| \leq |a_n-x_n|=r_n.$$
By Theorem \ref{T:local_symm}, the $T_{a_n}(S)-$ symmetric point
$$x_n^*=\sigma_{a_n}(x_n) \in supp f.$$
Using compactness of $supp f,$ choose a convergent subsequence
$$x_{n_k} \to x_0 \in supp f, \ k \to \infty.$$
Taking, if necessarily, a subsequence one more time, we can assume that also
$$r_{n_k} \to r_0.$$
Then, taking limits $a_n \to a, x_n \to x_0, r_n \to r_0,$ we will have
$$|a-x_0|=r_0$$
and for any $x \in supp f:$
$$|a-x| \leq r_0.$$
Those two inequalities show that
$$r_0=r,$$
where $r$ is defined in the formulation, and
$$x_0 \in K=K_{a,r}.$$
Now,
$$x_n^* =x_n-2 \langle x_n-a_n,\nu_{a_n} \rangle  \nu_{a_n} \to x_0-2 \langle x_0-a,\nu_{a} \rangle \nu_a =x_0^*,$$
as $ n \to \infty.$ Since $x_n^* \in supp f$ then $x_0^* \in supp f.$ Therefore $K^* \cap supp f \neq
\emptyset.$ Theorem is proved.

Theorems \ref{T:local_symm} and \ref{T:local_symm_mod} can be viewed as  non-linear
versions of the following global symmetry property, which follows from the uniqueness for
Cauchy problem for the wave equation:
\begin{Theorem} \label{T:CH} (\cite{CH}, Ch.VI, 8.1) Let $\Pi$ be a hyperplane in
$\mathbb R^d$ and $f \in C(\mathbb R^d).$ Then $Rf\vert_{\Pi \times (0,\infty)}=0$ if and
only if $f$ is odd with respect to reflections around $\Pi.$
\end{Theorem}
Obviously, $supp f$ in Theorem \ref{T:CH} is $\Pi-$ symmetric. Theorem
\ref{T:local_symm} states that if the hyperplane $\Pi$ is replaced by a hypersurface $S$
then, still, certain symmetry of $supp f$ holds, though in a much weaker (local) sense.

 The proof of Theorem \ref{T:antipodal} is geometric and is given in \cite{AQ1}. We
present it here to make the text of this article more self-sufficient.

\noindent
{\bf Proof of Theorem \ref{T:antipodal}}

We will present an analytic exposition of the geometric proof given in \cite{AQ1}.
We want to prove that if $f \in C_{comp}(\mathbb R^d)$ and $Rf(x,r)=0$ for all $x \in S$
and $r>0$ then $f=0$ or, equivalently, $supp f =\emptyset.$ We assume that $f \neq 0$ and
will arrive at a contradiction.

Since the tangent planes at $a$ and $b$ are parallel, the unit
normal vectors $\nu_a$ and
$\nu_b$ can be chosen equal
$$\nu_a=\nu_b=\nu=\frac{b-a}{|b-a|}.$$
Denote as above
$$\sigma_a(x)=x-2\langle x-a, \nu \rangle  \nu=x-2\frac{ \langle x-a,b-a \rangle }{|b-a|^2}(b-a). $$
the reflection around the tangent plane $T_a(S)$ and let $\sigma_b$ be the analogous
reflection for the point $b.$

Denote
$$r_1=max\{|x-a|: x \in supp f\}.$$

 Consider two cases:
\begin{enumerate}
\item $r_1 <|a-b|,$
\item $r_1 \geq |a-b|.$
\end{enumerate}

In the first case, $supp f$ lies on one side of $T_b(S):$
$$ \langle x-b,\nu \rangle < 0, \ x \in supp f$$
and therefore the entire $T_b(S)-$ symmetric set $\sigma_b(supp f)$  is disjoint from $supp f.$  This contradicts to Theorem \ref{T:local_symm_mod}.

Consider now the case $r_1 \geq |a-b|$ and  denote
$$r_2:=\sqrt{r_1^2-|a-b|^2}.$$
We claim that $supp f \subset \overline {B(b,r_2)},$ i.e. $|x-b| \leq r_2$ for all $x \in
supp f.$ To rpove that, consider
$$r=max\{|x-b|: x \in supp f\}.$$
Then $supp f \subset B(b,r)$ and it suffice to prove that $r \leq r_2.$

Suppose that $r > r_2.$  Denote
$$K=K_{b,r}=supp f \cap \{x \in \mathbb R^d: |x-b|=r \}.$$

 By Theorem \ref{T:local_symm_mod}, $K^*=\sigma_b(K)$ meets
$supp f.$  That means that there is $x_0 \in K$ such that $\sigma_b(x_0) \in K,$
i.e.,
$$x_0 \in supp f, \ |x_0-b|=r \ \mbox{and} \ x_0^*=\sigma_b(x) \in supp f.$$

Since $x_0 \in supp f$ then by definition of $r_1:$
$$|x_0-a| \leq r_1.$$
Therefore,
$$r_1^2 \geq |x_0-a|^2=\langle x_0-b+(b-a),x_0-b+(b-a) \rangle =|x_0-b|^2+|b-a|^2+2 \langle x_0-b,b-a \rangle.$$

Taking into  account that $$|x_0-b|=r, |b-a|^2=r_1^2-r_2^2,$$
we obtain the inequality
$$\langle x_0-b,b-a \rangle = r_2^2-r^2  <0.$$
But the same applies to the symmetric point $x_0^*=\sigma_b(x_0)$ because $x_0^*$ meets the same conditions $x_0^* \in supp
f$ and $|x_0^*-b|=|x_0-b|=r.$  Thus, also
$$\langle x_0^*-b,b-a \rangle <0.$$
Substitution
$$x_0=\sigma_b(x_0)=x_0-2\frac{\langle x-b,b-a \rangle }{|b-a|^2}(b-a)$$
yields
$$- \langle x_0 -b, b-a \rangle <0.$$
The obtained contradictions shows that $r \leq r_2$ and hence
$$supp f  \subset \overline B(b,r) \subset \overline{B(b,r_2)}.$$
Then we  repeat the argument, replacing $a$ by $b$ and $r_1$ by $r_2,$ and obtain
$$supp f \subset \overline{B(a, r_3)},$$
where $r_3=\sqrt{r_2^2-|a-b|^2}.$

Proceeding this way, we construct the sequence
$$r_{n+1}=\sqrt{r_n^2-|a-b|^2}, \  \mbox{i.e., } \ r_n=\sqrt{r_1^2-(n-1)|a-b|^2}, $$
such that
$$supp f \subset \overline{B(a,r_{2k+1})}, \ supp f \subset \overline{B(b,r_{2k})}.$$
For large enough $n$ we will have $r_n <|b-a|$ which, as it explained above, is
impossible. Therefore, the only possible conclusion is that $supp f=\emptyset$ and $f=0.$
Therefore, $S$ is an injectivity set.

\section{ Ruled surfaces}

Let $S$ be a real analytically ruled surface in $\mathbb R^3$ (see Definition
\ref{D:ruled_surface}). In accordance with the definition, $S$ consists of straight
lines, intersecting the fixed base curve $\gamma.$

More precisely, $S$ is locally the image of a map
$$ (t,\lambda) \to u(t,\lambda)= u(t)+\lambda e(t),$$
where
$$u(t): I  \to \mathbb R^3, \ e(t): I  \to S^2, \  I=(-1,1),$$
are real analytic vector-functions.

We denote $L_t$ the straight line
$$L_t=\{u(t)+\lambda e(t), \ \lambda \in \mathbb R\}.$$
\begin{Lemma} \label{L:ort} The parametrizing mapping $u(t)$ of the base curve $\gamma$
can be chosen so that
the tangent vector to the base curve and the directional vector are orthogonal:
\begin{equation}\label{E:ortog}
 \langle u^{\prime}(t), e(t) \rangle =0, \ t \in (-1,1).
\end{equation}
\end{Lemma}

\noindent
{\bf Proof} For any function $\lambda(t)$ we have
$$u(t,\lambda)=u(t)+\lambda(t)e(t) +(\lambda -\lambda(t))e(t).$$
Then $\mu=\lambda-\lambda(t)$ is a new parameter on the line $u(t)+\mathbb R e(t)$ and
therefore $S$ is the image of the mapping $\widehat u(t,\mu)= \widehat u(t)+\mu e(t),$
where $\widehat u(t)= u(t)+\lambda(t)e(t).$

The function $\lambda(t)$ is to be found from the condition
$$\langle \widehat u(t)^{\prime}, e(t)\rangle =
\langle u^{\prime}(t)+\lambda^{\prime}(t)e(t)+\lambda(t) e^{\prime}(t), e(t)\rangle
=\langle u^{\prime}(t),e(t) \rangle +\lambda^{\prime}(t)=0.$$ We have used here the  that
$\langle e(t),e(t) \rangle =1$ and $\langle e^{\prime}(t), e(t) \rangle =0.$ Therefore $\lambda(t)$ can be taken
$$\lambda(t)=-\int\limits_{t_0}^{t} \langle u^{\prime}(t),e(t) \rangle dt.$$
The condition of real analyticity preserves for $u(t)+\lambda(t)e(t).$

From now on, we assume that the parametrization $u(t,\lambda)$ satisfies  the
orthogonality condition (\ref{E:ortog}).

\subsection{ Regularity of the line foliation at smooth points}

In this subsection, we will prove that the line foliation of $S$ is regular at the points
where the surface $S$ is differentiable.

Notice that, in Definition \ref{D:ruled_surface}, the parametrizing mapping
$u(t,\lambda)$ is not assumed necessarily regular, i.e. the condition nondegeneracy of
the Jacobi matrix may be not fulfilled.

\begin{Definition} \label{D:regul} We call a point $a \in S$ of a ruled surface $S \mathbb R^3$  {\bf regular with respect to  a parametrization}
$ I \times I:  \ni (s,\sigma) \to w(s,\sigma), a=w(0,0),$  where $I=(-1,1),$ if
\begin{enumerate}
\item The mappings $ \mathbb R \ni \sigma \to w(s,\sigma)$ parametrize the same line foliation of $S.$
\item The mapping $w(s,\sigma)$ is differentiable and regular at $(0,0),$  i.e., the partial derivatives $\partial_s w(0,0), \partial_{\sigma}w(0,0)$ are linearly independent ( and then span the tangent space $T_a(S)$).
\end{enumerate}
We will call $a$ just {\bf regular point} of the given line foliation, if $a$ is regular with respect to some parametrization $w(s,\sigma).$

\end{Definition}.

\begin{Lemma} \label{L:smooth_ruling}
Let  $S_0$ be a ruled
surface with  $C^1 $ open base curve $W \subset S_0,$ i.e., $S_0=\cup_{w \in W}
L_w,$ where $L_w$ is a straight line passing through the point $w \in W.$ Suppose that $L_w \perp W, w
\in W.$ If $S$  is a$C^1-$ near a point $a \in W$ then $a$ is a regular point of the the
foliation $\{L_w, w \in W\}.$
\end{Lemma}

{\bf Proof}
 Let $\Omega_a$ be the neighborhoood of $a$ where $S$ is $C^1.$
$$I \ni s \to w(s) \in W$$ where $I$ is an open interval, be a $C^1$ parametrization of the base curve $W,$ and $\tau(w(s))=w^{\prime}(s)-$ the tangent vector to $W.$

Let $\nu(x), x \in \Omega_a,$ be the unit normal $C^1$ vector field on $\Omega_a.$  The surface $S$ is differentiable at $a,$ hence the normal unit vector
$\nu(a)$ is well defined, and $\nu(x)$ is $C^1$ mapping on $\Omega_a.$

Then the
cross-product
$$E(w)=\nu(w) \times \tau(w)$$
is both orthogonal to $W$ and tangent to $S_0$ and hence $E(w)$ is the directional vector of the
generating line $L_w.$ The vector field $E(w), w \in W$ is $C^1.$ Let
$$I \ni s \to w(s) \in W$$ where $I$ is an open interval, be a $C^1$ parametrization of the base curve $W.$
Then the mapping
$$I \times I  \ni (s,\sigma) \to w(s,\sigma)=w(s)+\sigma E(s), \sigma \in \mathbb R^3,$$ where
$$e(s):=e(w(s)),$$ parametrizes the given line foliation $\{L_w\}$ and satisfies  Definition \ref{D:regul} of  regular point.

Indeed, $w(s,\sigma)$ is differentiable at $(0,0),$ because $w(s)$ and $E(w(s))$ are differentiable.  The vectors
$$\partial_s w(0,0)=\tau(w), \  \partial_{\sigma}w(0,0)=E(0)$$
are nonzero and orthogonal to each other, hence the point $(0,0)$ is regular with respect to the parametrization $w(s,\sigma)$ of the given foliation,

Lemma is proved.

\section{ The structure of real analytically ruled  algebraic surfaces near singular points}

 In this section we study singular points of algebraic real-analytically ruled surfaces
in $\mathbb R^3.$  We did not find a relevant result in the literature, for , to our knowlegde, singular points of ruled surfaces and caustics of normal fields
(cf. \cite{AGV}) are classified for either generic surfaces or in the case of stable singularities, while in our situation, the surface and a point are given
and cannot be perturbed.

The following theorem, combined with Theorem \ref{T:antipodal}, will
be one of key points in the proof of the main result of this article:


\begin{Theorem} \label{T:conical} Let $(-1,1) \ni t \to u(t) \in \mathbb R^3$ and $(-1,1) \ni t \to e(t) \in S^2$
be two real analytic mappings. Denote $S$ the ruled surface $S:=\{u(t)+\lambda e(t), t
\in (-1,1), \lambda \in \mathbb R \}$ and assume that $S$ is algebraic.
Then the following four cases are possible:
\begin{enumerate}
\item $S$ is $C^1-$ manifold and the line foliation $\{L_t\}$ is regular at any point $a \in S.$
\item $S$ is a plane.
\item $S$ is a cone, i.e. all the lines $L_t$ have a common point (vertex).
\item $S$ has a {\bf cuspidal (double tangency) point} $a \in S$, which means the following:
if $H$ is a polynomial vanishing on $S$ and $H(x+a)=H_{k}(x)+H_{k+1}+...+H_N(x),$ where $H_j$ are homogeneous polynomials of degree $j$ and $H_{k} \neq 0 ,$ then  $H_k$ is divisible by a nonzero degenerate quadratic form $Q(x)=(A_1 x_1+A_2x_2+A_3x_3)^2.$
\end{enumerate}
\end{Theorem}

\begin{Remark} In both cases 1 and 2 $S$ is a smooth manifold, but in case 2, when  $S$ is a plane, the given line foliation can be  singular
(have caustics). For example, all the lines $L_t$ can pass through the same point, so that $S$ belongs to case 3, or there can be caustics of more complicated forms. On the other hand, planes can be viewed also as a regular ruled surface ( foliated into parallel lines) but this foliation can be different from the initial one .
\end{Remark}

\begin{Example} It was proved in \cite{IT} that  generic ruled
surface in $\mathbb R^3$ is equivalent, near its singular point, to Whitney umbrella,
which is the image $S$ of the mapping
$$(t,\lambda) \to (t^2,\lambda , \lambda t).$$
Whitney umbrella is algebraic surface with the equation $$z^2- yx^2=0.$$

The origin
$a=(0,0,0)$ is the only singular point.
Whitney umbrella is a typical ruled surface with cuspidal singular point, as  defined in case 4 of  Theorem \ref{T:conical}. Indeed, any
polynomial $H$ vanishing on $S$ is divisible by $x_3^2-x_1 x_2^2.$ Then the minor homogeneous
term $H_k$ of $H$ is divisible by $x_3^2,$ i.e., the property 3 holds with $Q(x_1,x_2,x_3)=x_3.$
\end{Example}

The important corollary of Theorem \ref{T:conical} is:
\begin{Corollary}\label{C:harm}
Let $S$ be as in Theorem \ref{T:conical}. Suppose that $S \subset H^{-1}(0),$ where $H$ is a nonzero harominic polynomial. Then $S$ is the surface of one of the first three cases in Theroem \ref{T:conical}.
\end{Corollary}

{\bf Proof} Suppose that $S$ is a surface of the four type, i.e, $S$ has  a cuspidal point $a \in S.$
Let $H$ be a harmonic polynomial such that $H \vert S=0.$ Then the minor term $H_k$ in the homogeneous decomposition
$$H(x+a)=H_k(x)+...HM(x)$$
is divisible by a nonzero quadratic polynomial $A^2(x)$ where $A(x)=A_1x_1+A_2x_2+A_3x_3$ is a nonzero linear form. Then
$$H_k(x)=0, \ \nabla H_k(x)=0, \mbox{whenever} \  A(x)=0. $$
Thus,  $H_k$ satisfies on the plane $\Pi=\{A(x)=0\}$ both the zero Dirichlet and Neumann
conditions. Since $H_k$ is harmonic,this implies $H_k=0$ identically. Therefore, the
homogeneous decomposition of $H$ begins with $H_{k+1}.$ The same argument yields
$H_{k+1}=0.$ Proceeding this way, we obtain $H=0.$  This contradiction shows that case 4 is impossible.

\subsection{Outline of the proof of Theorem \ref{T:conical}}

First of all, we will show that if  $a$ is not a conical point of $S$ then by a
suitable changing parameters $t$ (reparametrization) and $\lambda$ (rescaling), we can
pass to a parametrization (\ref{E:canonical}) of $S$ of the form
$$u(s,\sigma)= s^m v_m+
\sigma s^m e_0 + D(s,\sigma)\tau,$$ where $v_m, e_0, \tau$ are nonzero pairwise orthogonal
vectors and $D(s,\sigma)$ is a nonzero (if $S$ is not a plane) real analytic function.

Then we show that if $m$ is odd then $S$ is $C^1-$ differentiable at $a$ and, even more,
$a$ is a regular point of the line foliation on $S$ (Lemmas \ref{L:smo} and
\ref{L:smooth_ruling}).

In the case of even $m$ we reduce the situation,  by consequent descending the power $m$,
to the case of even $m$ and $D$ not even function of $s$ (we assume that $D \neq 0$ identically since otherwise $S$ is a plane).

Then we prove in Lemma
\ref{L:Deven} that in this case the point $a$ is of cuspidal type, i.e., the fourth case
of Theorem \ref{T:conical} takes place.

Thus, we conclude that if $S$ contains no cuspidal points then either $S$ is a plane or a cone, or the power $m$ associated with any point $a \in S$ is odd
and therefore $S$ is everywhere $C^1$ differentiable and the line foliation is everywhere regular.

\subsection{Preliminary constructions}\label{S:preliminary}

Let $a$ be a singular point of the real analytically ruled surface $S.$

As it is showed in Lemma \ref{L:ort}, we can choose the parametrization $u(t,\lambda)=u(t)+\lambda e(t)$ near $a$  so that
$\langle u^{\prime}(t),e(t) \rangle =0.$  Using translation we can always move $a$ to the origin and assume that $a=0.$ We can also assume that the value of
the parameter corresponding to the point $a$ is $t=0.$

\begin{Lemma}\label{L:sing} Let $a=u(0)+\lambda_0 e(0)=0$ be a singular point of the ruled surface $S.$ Then the parametrizing mapping $u(t,\lambda)=u(t)+\lambda t$ can be rewritten as $u(t,\mu)=u(t)+\mu e(t),$ where
\begin{equation}\label{E:v(t)}
v(t)=u(t)+\lambda_0 e(t), \ \mu=\lambda-\lambda_0
\end{equation}
and
\begin{enumerate}
\item
$v^{\prime}(0)=0.$
\item
If $v(t)=0$ identically then $S$ is a cone with the vertex $0.$ Otherwise, $v(t)$
decomposes in a neighborhood of $t=0$ into power series:
$$v(t)=v_m t^m+v_{m+1}t^{m+1}+..., \ v_m \neq 0,$$
where $m \geq 2,$  $v_j$ are vectors in $\mathbb R^3.$
\item $ \langle v_m,e(0) \rangle =0.$
\end{enumerate}
\end{Lemma}

{\bf Proof} Since $a$ is singular, the vectors
$$\frac{\partial u}{\partial t}(0,\lambda_0)=u^{\prime}(0)+\lambda_0 e^{\prime}(0),
\frac{\partial u}{\partial \lambda}(0,\lambda_0)=e(0)$$ are linearly dependent at
$0,\lambda_0:$

$$c_1 (u^{\prime}(0)+\lambda_0 e^{\prime}(0))+c_2 e(0)=0,$$
for some $c_1,c_2 \in \mathbb R, c_1^2+c_2^2 \neq 0.$

The unit vector $e(0)$ is orthogonal both to $u^{\prime}(0)$ and $e^{\prime}(0)$ ,
therefore $c_2=0$ and
$$u^{\prime}(0)+\lambda_0 e^{\prime}(0))=0.$$

Now rewrite $u(t,\lambda)$ as
$$u(t,\lambda)=u(t)+\lambda_0 e(t)+ (\lambda-\lambda_0)e(t),$$
and denote $\lambda-\lambda_0=\mu.$ Then we get the parametrization
$$u(t,\mu)=v(t)+\mu e(t), u(0,0)=a,$$
where
$$v(t)=u(t)+\lambda_0 e(t),$$
Then
$$v(0)=u(0)+\lambda_0 e(0)=0, \ v^{\prime}(0)=0.$$

The two cases are possible

1) $v(t) \equiv 0.$

Then $u(t,\lambda_0)=u(t)+\lambda_0 e(t)=v(t)=0,$ i.e., all the lines $L_t$ pass through
the origin and therefore $S$ is a cone with the vertex $0.$

2) $v(t)$ is not identical zero .

Then by real analyticity:
\begin{equation} \label{E:m}
u(t,\mu)= v_mt^m +.... + \mu (e_0 +e_1 t+....),
\end{equation}
where $ v_m \neq 0.$  Since $v^{\prime}(0)=0$ then  $m \geq 2.$

Also we have
$$\langle v^{\prime}(t), e(t) \rangle =\langle u^{\prime}(t)+\lambda_0 e^{\prime}(t),e(t) \rangle =0.$$
Thus,
$$\langle mv_mt^{m-1}+..., e_0+e_1t+... \rangle =0$$
and dividing by $t^{m-1}$ and letting $ t \to 0$ yields
$$ \langle v_m,e_0 \rangle =0.$$  Lemma is proved.

On the next step, we will replace the parameters $\mu,t$ by  new parameters $\sigma,s$ which are
more convenient for further investigation. We start with  re-scaling the parameter $\mu$
on the ruling lines.

\subsection{Re-scaling: changing the linear parameter $\mu.$}

Thus, by Lemma \ref{L:sing}, the surface $S$ is parametrized, near $a=0,$ by the mapping
$u(t,\mu)=v(t)+ \mu e(t),$ where
$$v(t)=\sum\limits_{j=m}^{\infty}v_j t^j, \ e(t)=\sum\limits_{j=0}^{\infty}e_j t^j$$
and $ \langle v_m,e(0) \rangle =0.$

Let $\tau$ be a unit vector orthogonal both to $v_m$ and $e_0.$ Then the triple
$$v_m,e_0,\tau$$
constitutes a basis in $\mathbb R^3.$

Decompose the vector-coefficients $v_m, v_{m+1},...$ and $e_0, e_1,..., $ into  linear
combinations of the basis vectors:
$$v_j=A_j v_m+B_j e_0+C_j \tau, \ j \geq m,$$
$$e_j=\hat A_j v_m+ \hat B_j e_0+ \hat C_j \tau, \ j \geq 0,$$
and since $v_m, e_0, \tau$ constitute the basis, one has
$$ A_m=1, B_m=0, C_m=0, \ \hat A_0=0, \hat B_0=1, \hat C_0=0.$$
Substitution the expressions for $v_j, e_j$ into the the power series for $v(t)$ and $e(t)$ leads to:
$$v(t)=A(t)v_m+B(t)e_0+C(t)\tau,$$
$$e(t)=\hat A(t)v_m+ \hat B(t)e_0+ \hat C(t)\tau,$$
where we have denoted
\begin{equation} \label{E:ABC}
A(t)=\sum\limits_{j=m}^{\infty} A_jt^j,  \ \ B(t)=\sum\limits_{j=m+1}^{\infty}B_j t^j, \ \ C(t)=\sum\limits_{j=m+1}^{\infty}C_j t^j
\end{equation}

and
\begin{equation} \label{E:hatABC}
\hat A(t)=\sum\limits_{j=1}^{\infty} \hat A_j t^j, \ \ \hat B(t)=\sum\limits_{j=0}^{\infty} \hat B_j t^j, \ \ \hat C(t)=\sum\limits_{j=1}^{\infty} \hat C_j t^j.
\end{equation}

Correspondingly, the  parametrizing function $u(t,\mu)=v(t)+\mu e(t)$  takes the form
\begin{equation}\label{E:parametrization}
u(t,\mu)=(A(t)+\mu \hat A(t))v_m + (B(t)+\mu \hat B (t))e_0+(C(t)+\mu \hat C(t))\tau.
\end{equation}

 Let us fix a real number $\sigma \in \mathbb R$ and write the functional equation
\begin{equation} \label{E:B}
B(t)+\mu \hat B(t)=\sigma (A(t)+\mu \hat A(t)).
\end{equation}
This equation defines the parameter $\mu$ as a function of $\sigma$ and $t:$
$$\mu=\mu(\sigma,t)=\frac{\sigma A-B}{\hat B- \sigma \hat A}.$$

Since from (\ref{E:ABC})
$$B(t) = B_{m+1}t^{m+1}+...;  \ \ \ \hat B(t) =1+B_1 t+... $$
and
$$ A(t)=t^{m}+A_{m+1} t^{m+1}+...; \ \ \ \hat A(t) =A_1 t+..., $$
and $m >1,$ we obtain
$$\mu=\frac{\sigma A-B}{\hat B- \sigma \hat A}=
\frac{\sigma t^m+...-B_{m+1}t^{m+1}+...}{(1+ \hat B_1 t+...)-\sigma ( \hat A_1t+...)}$$
and hence
$$\mu=\mu(t)=\sigma  t^m+ o(t^m).$$

Then the coefficient $A(t)+\mu(t,\sigma) \hat A(t)$ in front of $v_m$ in
(\ref{E:parametrization}) is

$$ A(t)+\mu(t,\sigma) \hat A(t)=A_m t^m+...+(A_m\sigma t^m+...)(A_1t+...)=A_m t^m+o(t), \ t \to 0.$$

\begin{Remark} \label {T:base} The base curve $\{t \to v(t)\}$ of the foliation is given by the condition $\mu=0$ which corresponds, due to (\ref{E:B}), to
$$\sigma=\frac{B(t)}{A(t)}=B_{m+1}t+o(t).$$
\end{Remark}

\subsection{ Re-parametrization: changing the parameter $t$ of the base curve.}

Now introduce  the new parameter $s$ by the relation
$$s^m=A(t)+\mu \hat A(t)=t^m+o(t), \ t \to 0.$$
If $m$ is odd, then the real parameter $s=s(t)$ is well defined near $t=0.$ If $m$ is
even then $s=s(t)$ near $t=0$ is the real branch of $(A(t)+\mu \hat A(t))^{\frac{1}{m}}$
for which
$$s=s(t)=t+o(t).$$
Thus, that asymptotic holds for both odd and even $m.$

From (\ref{E:parametrization}) and (\ref{E:B}), one can rewrite, in a neighborhood of $s=0,$ the function $u(t,\mu)$ as a function of the new
parameters $s, \sigma:$
\begin{equation}\label{E:canonical}
u(s,\sigma)= s^m v_m+ \sigma s^m e_0 + D(s,\sigma)\tau,
\end{equation}
where we have denoted
$$D(s,\sigma):=C(t)+\mu \hat C(t).$$

Since $s=t+o(t),$ we have from (\ref{E:ABC}),(\ref{E:hatABC}):
$$C(t)=C_{m+1}t^{m+1}+o(t^{m+1})=C_{m+1}s^{m+1}+o(s^{m+1}),$$
$$\hat C(t)=\hat C_1 t+o(t)=\hat C_1 s +o(s),$$
$$\mu=\sigma  s^m+o(s^{m}).$$

Then we have
\begin{equation}\label{E:D}
D(s,\sigma) =C(t)+\mu \hat C(t)= (C_{m+1}+\sigma \hat C_1)s^{m+1}+o(s^{m+1}).
\end{equation}

\begin{Lemma} \label{L:H_k=0} If $H$ is a polynomial vanishing on $S$ and
$H=H_k+H_{k+1}+..$ is its decomposition into homogeneous polynomials, then $H_k(x)=0$ for
all vectors $x \in span \{v_m, e_0 \}.$
\end{Lemma}

\noindent
{\bf Proof} We have $H(u(s,\sigma))=0$ for all $\sigma \in \mathbb R$ and $s$ close to 0.
From (\ref{E:D}),  $D(s,\sigma)=o(s^m)$ and then formula (\ref{E:canonical}) implies
$$H(u(s,\sigma))=H_k (s^m v_m+\sigma s^m e_0 + o(s^m))+H_{k+1}(s^m v_m+\sigma s^m e_0 + o(s^m))+...=0.$$
Since $H_j$ are homogeneous of degree $j$, dividing by $s^{km}$ and letting $s \to 0$
yields:
$$H_k(v_m +\sigma e_0)=0.$$
Then
$$H_k(\alpha v_m +\alpha \sigma e_0)=\alpha^k H_k(v_m+\sigma e_0)=0$$
for any $\alpha \in \mathbb R.$ Since $\sigma$ is arbitrary, the real numbers $\alpha, \alpha \sigma$ are arbitrary as well, and hence $H$ vanishes on any linear
combination of the vectors $v_m$ and $e_0.$ Lemma is proved.

\begin{Lemma}\label{L:D=0} If $D(s,\sigma)=0$ identically then $S$ locally is a plane ( case 1 of
Theorem \ref{T:conical}).
\end{Lemma}

\noindent
{\bf Proof} If $D(s,\sigma) \equiv 0$  then we have from (\ref{E:canonical})
$u(s,\sigma)=s^m v_m+\sigma e_0$ and hence the image of $u$ is contained in the plane
spanned by the vectors $v_m$ and $e_0.$

\begin{Lemma} \label{L:reduce} Suppose that $D(s,\sigma)$ is not identically zero. Then a suitable  change of the parameter $s$ leads to
one of the following cases: hold for the power $m$ in (\ref{E:canonical}) :
\begin{enumerate}
\item The integer $m$ in (\ref{E:canonical}) is odd.
\item $m$ is even but $D(s, \sigma)$ is not an even function with respect to $s.$
\end{enumerate}
\end{Lemma}

\noindent
{\bf Proof} We will consequently descend the power $m$ until we reach one of the above cases.

If $m$  is odd then we are done. Suppose that
$m$ is even, $m=2m^{\prime}.$ If $D(s,\sigma))$ is an not even with respect to $s$, then
we are done.

If $D(s,\sigma)$ is still even in $s$ then $D(s)=D^{\prime}(s^2),$ where $D^{\prime}$ is
a new function, real analytic in $s$ near 0.

Then introduce new parameter
$$s^{\prime}=s^2$$
and pass to the new parameter $s^{\prime}$ and the new parametrizing function
$$u(s^{\prime},\sigma)=(s^{\prime})^{m^{\prime}}v_m +(s^{\prime})^{m^{\prime}}
e_0+D^{\prime}(s^{\prime})\tau,$$ which extends as a real analytic function to negative
values of $s^{\prime}.$

If,  again, both functions $(s^{\prime})^{m^{\prime}}$ and $D^{\prime}(s^{\prime})$ are
even, we introduce the new parameter
$$s^{\prime\prime}=(s^{\prime})^2.$$

Proceeding that way, we finally end up either with odd $m$ or with even $m$  but not even
(with respect to $s$)  function $D(s,\sigma).$ Lemma is proved.

\subsection{The case of even $m$}

The following lemma shows that the case of even power $m$ leads to the case 4 in Theorem
\ref{T:conical}, of double tangency at the singular point $a$ (which here is assumed to be
$a=0$):

\begin{Lemma} \label{L:Deven} Let $m$ be even and let $D(s,\sigma)$ be not identically zero function
(i.e. due to (\ref{E:canonical}) the surface $S$ is not a plane). Then $a$ is a cuspidal point as defined in case 4 of Theorem \ref{T:conical}/
\end{Lemma}

\noindent
{\bf Proof}

\noindent

\subsubsection{ Extracting the even part of $D(s,\sigma)$}

 By Lemma \ref{L:reduce} we can make, by means of a suitable reparametrization,  the
function $D(s,\sigma)$ not even  with respect to the variable $s.$

We fix an arbitrary $\sigma$ such that $D(s,\sigma)$ is not even in $s.$ By the
construction, the power series for $D$ contains no powers of $s$ less than $m+1:$
$$D(s,\sigma)=\sum\limits_{j=m+1}^{\infty}D_j(\sigma)s^j.$$
Since $D$ is not an even function with respect to $s,$ there exists at least one odd
exponent $j$ with $D_j(\sigma) \neq 0$ near $\sigma=0.$
Denote
$$j_0=min \{j \geq m+1: j \ \mbox{is odd and } \ D_j(\sigma) \neq 0 \}.$$

Let us split the above power series into two parts:
$$D(s,\sigma)=D_1(s,\sigma)+D_2(s,\sigma),$$
where
$$D_1(s,\sigma)=\sum\limits_{j=m+1}^{j_0-1} D_{1,j}(\sigma)s^j,$$
$$D_2(s,\sigma)=\sum\limits_{j=j_0}^{\infty} D_{2,j}(\sigma)s^j.$$

Then $D_1$ is even:
$$D_1(-s,\sigma)=D_1(s,\sigma)$$
because all the powers $j=m+1,...,j_0-1$ are even.

Now, we have:
$$D_1(s,\sigma)=D_{1,m+1}s^{m+1}+o(s^{m+1})$$
and
$$ D_2(s,\sigma)=D_{2,j_0}s^{j_0}+o(s^{j_0}),$$
with $j_0$ odd. It is important that
\begin{equation} \label{E:D_2}
D_{2,j_0} \neq 0.
\end{equation}

Substituting the above representations for $D_1(s,\sigma)$ and $D_2(s,\sigma)$ into
formula (\ref{E:canonical}) for $u(s,\sigma)$ we obtain
\begin{equation}\label{E:u}
u(s,\sigma)=s^m v_m+ \sigma s^m e_0+ (D_1(s,\sigma)+D_2(s,\sigma))\tau.
\end{equation}

\subsubsection{ Taylor series for $H(u(s,\sigma))$}

Now, let $H$ be  a polynomial vanishing on $S:$
$$H(x)=0 \ \forall x \in S.$$
We want to prove that $S$ has a double tangency at $a,$  more precisely, that the
property 2) of Theorem \ref{T:conical} is satisfied for the polynomial $H.$

From the representation (\ref{E:u}) ,we have
$$H(u(s,\sigma))=H(s^m v_m+ \sigma s^m e_0+ (D_1(s,\sigma)+D_2(s,\sigma))\tau))=0.$$

Now, let us write Taylor formula for the polynomial $H,$ at the point
$$s^m v_m + \sigma s^m e_0+ D_1(s,\sigma) \tau,$$
on the vector
$$D_2(s,\sigma)\tau.$$
It yields:
\begin{equation} \label{E:1st}
H(u(s,\sigma))=\sum\limits_{r=0}^{deg H} d^r H(s^m v_m+ \sigma s^m e_0+
D_1(s,\sigma)\tau; D_2(s,\sigma)\tau)=0,
\end{equation}
where $d^rH(a;h)$
stands for the $r-$th differential of $H$  at a point $a,$ evaluated on a vector $h.$

Replacing $s$ by $-s,$  we have, taking into account that $D_1$ is even in $s,$ one more
relation:

\begin{equation} \label{E:2nd}
H(u(-s,\sigma))=\sum\limits_r d^rH(s^m v_m+ \sigma s^m e_0+ D_1(s,\sigma)\tau,
D_2(-s,\sigma)\tau)=0.
\end{equation}

 Now, if we subtract the second identity from the first one, then the term
corresponding to $r=0$ cancels and we will have:
\begin{equation}\label{E:subtract}
H(u(s,\sigma))-H(u(-s,\sigma))=\sum\limits_{r=1}^{deg H} T_r=0
\end{equation}
where we have denoted
\begin{equation}\label{E:subtract1}
\begin{split}
T_r &=d^rH(s^m v_m+ \sigma s^m e_0+ D_1(s,\sigma)\tau;  D_2(s,\sigma) \tau ) \\
    &- d^rH(s^m v_m+\sigma s^m e_0+ D_1(s,\sigma)\tau; D_2(-s,\sigma)\tau).
\end{split}
\end{equation}
Here we have used that $m$ is even and $D(-s,\sigma)=D(s,\sigma).$

\subsubsection{Contribution of the first differential}

 Now let us look at the first term $T_1$ in
the expression (\ref{E:subtract})- (\ref{E:subtract1}), corresponding to the first
differential of $H:$

\begin{equation}\label{E:T1}
T_1= \langle \nabla H(s^m v_m+ \sigma s^m e_0+ D_1(s,\sigma)\tau),
(D_2(s,\sigma)-D_2(-s,\sigma))\tau \rangle + \mbox{{\it higher order differentials}}.
\end{equation}
Notice that  the first term in the power series for $D_2(s,\sigma)$ is
$D_{2,j_0}(\sigma)s^{j_0},$ where $j_0$ is odd. Therefore,
\begin{equation}\label{E:D_two}
D_2(s,\sigma)-D_2(-s,\sigma)=2D_{2,j_0} s^{j_0}+o(s^{j_0}).
\end{equation}
Also,
\begin{equation}\label{E:D_one}
D_1(s,\sigma)=D_{1,m+1}s^{m+1}+o(s^{m+1}).
\end{equation}

Now decompose $H$
$$H=H_k+....+H_{deg H}$$
into sum of homogeneous polynomials, $deg H_j=j,$ and substitute the decomposition into
(\ref{E:T1}):
$$T_1:=dH_k(...)-dH_k(...) + dH_{k+1}(...) -dH_{k+1}(...) + \mbox {higher order differentials}.$$
Here all the differentials $dH_k$ are evaluated at the point
$$s^m v_m+ \sigma s^m e_0+D_1(s,\sigma)\tau$$
and on the vector
$$D_2(\pm s,\sigma) \tau,$$
depending whether we have $+$ or $-$ in front of $dH_k$ in (\ref{E:T1}).

Now using (\ref{E:D_two}), (\ref{E:D_one}) and homogeneity of $H_k$ we obtain
$$T_1= s^{(k-1)m+j_0} \langle \nabla H_k (v_m+\sigma e_0+(D_{1,\sigma}s+o(s))\tau,
(2D_{2,j_0}+o(s))\tau \rangle + s^{km+j_0} \langle \nabla H_{k+1}(...),... \rangle +... .$$
and at last
\begin{equation} \label {E:T}
T_1=2D_{2,j_0} s^{(k-1)m+j_0} \langle \nabla H_k (v_m+\sigma e_0), \tau) \rangle +o(s^{(k-1) m+j_0}).
\end{equation}
Similarly, substituting  the above asymptotic (\ref{E:D_two}),(\ref{E:D_one}) of $D_1$
and $D_2$ into the next  homogeneous terms $H_{k+1}, H_{k+2},...$ leads to the
expressions similar to (\ref{E:T}) were  $k$ is replaced by $k+1$, $k+2$ and so on.
Therefore, the least power that comes from $H_{k+1}, H_{k+2},...$ is $ s^{k m+j_0}.$

\subsubsection{ Contribution of the higher differentials}

 Let us turn now to the higher
differentials and consider the contribution of the terms corresponding to  $d^2 H_k, d^3
H_k...$ in the asymptotic near $s=0.$

Consider now the term $T_2$ in (\ref{E:subtract}), corresponding to the second
differential $d^2H:$

$$d^2 H(s^m v_m+ \sigma s^m e_0+ (D_{1,\sigma}s^{m+1}+o(s^{m+1}) \tau ; (2D_{2,j_0} s^{j_0}+
o(s^{j_0}))\tau)).$$

The asymptotic of (\ref{E:subtract}) near $s=0$ is determined again by the minimal degree
homogeneous polynomial $H_k,$ more precisely, by the difference

$$d^2 H_k(s^m v_m+ \sigma s^m e_0+ (D_{1,\sigma}s^{m+1}+o(s^{m+1}) \tau,(2D_{2,j_0}s^{j_0}+o(s^{j_0}))\tau)),$$

which comes from the minor homogeneous term $H_k$ in $H.$

By the homogeneity, it equals to
$$4D_{2j_0}(\sigma) ^2 s^{(k-2)m+2j_0}d^2 H_k( v_m+ \sigma e_0+ o(s),\tau))+o(s^{(k-2)m+2j_0}).$$

However,
$$(k-2)m+2j_0=(k-1)m+j_0-m+j_0>(k-1)m+j_0,$$
because $j_0-m>0.$

Moreover, for the next terms, coming from the higher differentials $d^r,$ we will have the
following order of the asymptotic
$$(k-r)m+r j_0=(k-1)m+ j_0 -(r-1)m+(r-1)j_0>(k-1)m+j_0.$$

Thus,  we see that only the first differential $d H_k$ of the minor homogeneous term
$H_k$ contributes the term  $s^{ (k-1)m+j_0}$ of the minimal power to the asymptotic of
$H(u(s,\sigma))$ near $s=0.$

Therefore, the main term of the asymptotic, which is determined by the minimal power of
$s$ , equals to
$$H(u(s,\sigma))-H(-s,(-s,\sigma))= 2D_{2,j_0}(\sigma) s^{(k-1)m+j_0} \langle \nabla H_k (v_m+\sigma e_0), \tau \rangle  +... $$

\subsubsection{ Double tangency property}

 Since the left hand side is identically zero
$$H(u(s,\sigma))-H(-s,u(-s,\sigma))=0,$$
the main term of the asymptotic is zero as well. It follows then from $D_{2,j_0} \neq 0$
that
$$ \langle \nabla H_k (v_m+\sigma e_0), \tau \rangle =0.$$

Now recall that $\sigma$ is an arbitrary real number. Since the polynomial $H_k$ is
homogeneous, we have
$$\langle \nabla H_k (h), \tau \rangle =0,$$
for all $h \in \Pi:=span \{v_m,e_0 \}.$ Since the vector $\tau$ is orthogonal to the
plane $\Pi,$ the normal derivative
$$\frac{\partial H_k}{\partial \tau}=0$$
on $\Pi.$

Also, we know from Lemma \ref{L:H_k=0} that $H_k=0$ on $\Pi.$ Thus $H_k$ vanishes on
$\Pi$ at least to the second order and therefore if to define linear form
$$A(x)=\langle (x, \tau \rangle , $$
then $H$ is divisible by $Q^2:$
$$H=A^2R.$$  Lemma is proved.

\subsection {The case of odd $m$}

\begin{Lemma}\label{L:smo} If $m$ is odd then  the surface $S$ is differentiable at $a=0.$ If, moreover,
$S$ is differentiable in a neighborhood of the point $a$ then $S$ is a $C^1-$ manifold there.

\end{Lemma}

\noindent
{\bf Proof}
By Lemma \ref{L:sing}, the surface $S$ is  the
image of the function
$$u(t,\mu)=v(t)+\mu e(t),$$
where
$$ v(t)= v_m t^m+... ; \ e(t)= e_0+ e_1 t+..., \ m=2s+1.$$
Since $m$ is odd, the
curve parametrized by
$$u(t,0)=v(t), t \in I=(-\varepsilon,\varepsilon),$$
is differentiable, which follows from the change of the parameter $t^m=s:$
$$ v(s)=v_m s +v_{m+1} s^{\frac{m+1}{m}}+...=v_m s + o(s).$$
We also have from definition (\ref{E:v(t)}) of $v(t)$ and Lemma \ref{L:ort}:
$$\langle v^{\prime}(t),e(t) \rangle =\langle u^{\prime}(t)+\lambda_0 e^{\prime}(t), e(t) \rangle =0.$$
Therefore, the image of the function $u(t,\mu)$ describes a ruled surface consisting of
straight lines orthogonal to the differentiable curve $v:I \to \mathbb R^3.$

Apply an orthogonal transformation so that the  triple $v_m, e_0, \tau$ becomes the axis.  Denote $x_1,x_2,x_3$ the coordinates of points in the basic
$v_m, e_0, \tau.$

Then, according to (\ref{E:canonical}), the mapping $u(s,\sigma)$ has the following
representation in the new coordinates:
$$u(s,\sigma)=(x_1,x_2,x_3)=(s^m, \sigma s^m, D(s,\sigma)).$$
We have

\begin{eqnarray*}
x_1 &=& s^m, \\
x_2 &=& \sigma s^m, \\
x_3 &=& D(s,\sigma),
\end{eqnarray*}

and therefore
$$s=x_1^{\frac{1}{m}}, \sigma=\frac{x_2}{x_1}.$$
The function $D(s,\sigma)$ is real analytic at $s=0, \sigma=0:$
$$D(s,\sigma)=\sum_{\alpha,\beta \in Z_{+}}c_{\alpha,\beta} s^{\alpha}\sigma^{\beta},$$
in a neighborhood of $ s=0, \sigma=0.$

Moreover, according to (\ref{E:D}), $D(s,\sigma)=o(s^m), s \to 0$ and hence
$$\alpha \geq m+1$$
in the Taylor series for $D.$

Substituting the expressions for $s,\sigma$ through $x_1,x_2$ yields the representation of the
function
$$x_3=z(x_1,x_2)=D(x^{\frac{1}{m}}, \frac{x_2}{x_1})$$ as a Newton-Puiseux fractional
power series:
\begin{equation}\label{E:series}
z(x_1,x_2)=\sum_{\alpha=m+1,\beta=0}^{\infty}
c_{\alpha,\beta}x_1^{\frac{\alpha}{m}-\beta}x_2^{\beta}.
\end{equation}

\subsubsection{ Differentiability  of $z(x_1,x_2)$ at $(0,0)$}

We know that the line $L_0=\{\lambda e_0, \lambda \in \mathbb R\}$ is one of the generating lines and belongs to $S.$
In the coordinates $x_1,x_2,x_3$, the line $L_0$ has the equation $x_1=x_3=0.$  Since $x_3=z(x_1,x_2)$ is the equation of $S,$
we conclude that
$$\lim\limits_{x_1 \to 0}z(x_1,x_2)=0$$
for any fixed $x_2.$  This implies that the series (\ref{E:series}) contains
only positive powers of $x_1.$

 Therefore,
the series can be rewritten as

\begin{equation}\label{E:Puis}
z(x,y)=\sum\limits_{\nu>0,\beta \geq 0} b_{\nu,\beta} x_1^{\nu}x_2^{\beta},
\end{equation}
where we have introduced the new coefficients
$$b_{\nu,\beta}=c_{\alpha,\beta}, \  \nu=\frac{\alpha}{m}-\beta.$$
In our case $\nu$ is strictly positive because $z(0,0,0)=0.$

Notice, that since $m$ is odd, the fractional power $x_1^{\nu}$ is well defined for $x_1<0$
as well, so the decomposition (\ref{E:Puis} ) holds in a full neighborhood of $(0,0).$

The general term in the Newton-Puiseux series (\ref{E:Puis})  is of homogeneity degree
$$\nu+\beta=(\frac{\alpha}{m}-\beta)+\beta=\frac{\alpha}{m} > 1+\frac{1}{m}.$$

The series (\ref{E:Puis}) can be written in the polar coordinates
$$x_1=r cos \theta, \ x_2=r \sin \theta $$ as
$$z(x_1,x_2)=\sum\limits_{nu > 0,\beta \geq 0}b_{\nu,\beta} r^{\nu+\beta} (cos \theta)^{\nu} (sin \theta)^{\beta}.$$

Since the exponents $\nu,\beta \geq 0$ then $|cos \theta|^{\nu}, |sin \theta|^{\beta} \leq 1,$ the inequality $$\nu+\beta >1 +\frac{1}{m}$$ implies
$$z(x_1,x_2)=o(r), \  r \to 0.$$
Therefore the function $z(x_1,x_2)$ is differentiable at $(0,0)$ with $dz(0,0)=0.$ Lemma is proved.

\subsubsection{ $S$ is differentiable in a neighborhood of $a$ implies $S$ is $C^1$}

 If $S,$ which is the graph of the function $z(x_1,x_2)$ is differentiable in  a neighborhood of $a=0$ then the $z(x_1,x_2)$ is differentiable at any point in a neighborhood $U$ of $(0,0).$  Due to (\ref{E:Puis})
 \begin{equation}\label{E:Puis1}
 \frac{\partial z}{\partial x_1} (x_1,x_2)=\sum\limits_{\nu>0,\beta \geq 0} b_{\nu,\beta} \nu x_1^{\nu-1}x_2^{\beta}.
 \end{equation}
 Since $\nu > 0 $ is fractional, the number $\nu-1$ can be negative. However, this is not the case, because
 if  series (\ref{E:Puis1}) contains negative powers of $x_1$ then for small $x_2 \neq 0$ we have
 $\lim\limits_{x_1 \to 0} \frac{\partial z}{\partial x_1} (x_1,x_2)=\infty$ which contradicts to the differentiability of $z(x_1,x_2)$ at the points $(0,x_2)$ with small $x_2.$
 Then the series
 $$\frac{\partial z}{\partial x_1} (x_1^m,x_2)=\sum\limits_{\nu>0,\beta \geq 0} b_{\nu,\beta} \nu x_1^{m(\nu-1)}x_2^{\beta}$$
 is a power series, since $m(\nu-1) =\alpha- m\beta-m$ is integer and nonnegative. Power series are continuous
 in their domains of convergence, therefore $\frac{\partial z}{\partial x_1} (x_1^m,x_2)$ is continuous in a neighborhood of $(0,0).$  Since $m$ is odd, the mapping $x_1 \mapsto x_1^m$ is a homeomorphisms and hence the continuity of $\frac{\partial z}{\partial x_1} (x_1^m,x_2)$ follows.

 Same argument implies the continuity of the $\frac{\partial z}{\partial x_2}$ since $x_2$ the series (\ref{E:Puis}) in just a usual power series with respect to $x_2.$  The proof of Lemma is completed.




\subsection{End of the proof of Theorem \ref{T:conical}}
Now we are ready to finish the proof of Theorem \ref{T:conical}.

We start with assumption that $S$ is neither a plane nor a cone. Then we have to prove that either the surface $S$ is $C^1$ manifold and the line foliation is everywhere regular or $S$ has a cuspidal point.

Lemma \ref{L:even} says that cuspidal singular points $a \in S$ correspond  to the even associated powers $m$ in decomposition (\ref{E:canonical}).  Therefore, if $S$ is free of cuspidal points, then for any singular (with respect to the initial parametrization of our line foliation) point the associated power is odd.

But Lemma \ref{L:smo} implies that then $S$ is differentiable at any singular  point $a \in S.$ Surely, $S$ is also differentiable at any regular point. Therefore $S$ is differentiable everywhere. But then the second assertion of Lemma \ref{L:smo} yields that $S$ is $C^1-$ manifold and the line foliation of $S$ is everywhere regular (with respect to some parametrization of the line foliation).

  Thus, we have proven that one of the fourth cases enlisted in Theorem \ref{T:conical}  holds. Theorem is proved.
\section {Irreducible case.Proof of Theorem \ref{T:Main1}}

\subsection{Extremal ruling lines and antipodal points}

 For any two ruling straight lines $L_t, L_s \subset S$ define the distance function
$$d(t,s):=dist(L_t,L_s)=min\{|u-v|:u \in L_t, v \in L_s \}.$$
\begin{Lemma} If $d(s,t)=0$ for all $t,s$ then  $S$ is a cone.
\end{Lemma}

\noindent
{\bf Proof} The condition implies that any two ruling lines meet. Fix two non-parallel
ruling lines $L_t, L_s.$ They intersect at some point $a \in L_t \cap L_s.$

Due to real analyticity of the one-dimensional connected family $\{L_t\}$ of the ruling
lines, the two cases are possible:

1) all  the lines $L_t$  pass through the point $a$, and then $S$ is a cone with the
vertex $a,$

2) at most finite number of lines $L_{t_1}, \cdots L_{t_N}$ contain $a.$

Suppose that case 2) takes place. Take any third ruling line $L_r$ for $r \neq
t_1,...,t_N.$ Since any two ruling lines have a common point, the line $L_r$ must
intersect both lines $L_t, \ L_s$ at points different from $a.$ This implies that $L_r$
belongs to the two-dimensional plane $\Pi$ spanned by $L_t, L_s.$ Therefore all but at
most finite number of ruling lines belong to $\Pi.$ This implies that the union of those
lines $S=\Pi.$ Therefore $S$ is a 2-plane which, of course, is a cone.

Now we are interested in the case when $d(s,t)$ is not identically zero function.
\begin{Lemma}\label{L:maximum}
If $S$ is not a cone then there are two maximally distant ruling lines $L_{t_0},
L_{s_0},$ i.e., the distance function $d(t,s)$ attains its maximum:
$$d(t_0,s_0)=\max \limits_{t,s}d(t,s) >0.$$
at some values $t_0,s_0$ of the parameters.
\end{Lemma}

\noindent
{\bf Proof} The function $d(s,t)$ is defined on the compact set $[-1,1] \times [-1,1].$
It is upper semi-continuous, i.e., the upper limit
$$\limsup d(t,s)_{(t,s)\to (t_0,s_0)} \leq
d(t_0,s_0).$$

Indeed, let $a=u(t_0)+\lambda_0 e(t_0) \in L_{t_0}, \ b=u(s_0)+\mu_0 e(s_0) \in L_{s_0},$
be the points on the straight lines $L_{t_0}, L_{s_0}$ such that
$$|a-b|=dist(L_{t_0},L_{s_0}).$$

If $(t_n,s_n) \to (t_0,s_0)$ then
$$a_n=u(t_n)+\lambda_0 e(t_n) \to a, b_n=u(s_n)+\mu_0 e(s_n) \to b.$$

Then we have $$d(s_n,t_n) \leq |a_n-b_n|$$ and hence
$$\lim\limits_{n \to \infty} d(t_n,s_n) \leq \lim\limits_{n \to \infty}
|a_n-b_n|=|a-b|=d(t_0,s_0).$$

Due to the arbitrariness of the sequence $(t_n,s_n) \to (t_0,s_0),$ the function $d(t,s)$
is upper semi-continuous. By Weierstrass theorem it attains its maximal value
$d(t_0,s_0)$. Since $d(t,s)$ is not identically zero function, we have $|a-b|=d(t_0,s_0)
>0.$ We will call $a,b$ {\bf extremal points.}

\begin{Lemma}\label{L:anti} Suppose that the line foliation of $S$  contains no parallel lines.
Suppose that the surface $S$ is differentiable at the extremal points $a$ and $b$ and the foliation $S=\cup_t L_t$ is regular at both extremal points $a$ and $b.$ Then $a$ and $b$ are antipodal points (see Definition \ref{D:antipodal_points}.
\end{Lemma}
According to Definition \ref{D:regul}, regularity means that near the points $a$ and $b,$ the surface $S$ is the image of
of the mappings
$$w_a(t,\lambda)=w_a(t)+\lambda E_a(t), \ w_b(s,\mu)=w_b(s)+\lambda E_b(s),$$
correspondingly, which define the same foliation and are differentiable and regular at the points $(t_0,\lambda_0), (s_0,\mu_0).$
Here $a=u_a(t_0,\lambda_0), \ b =u_b(s_0,\mu_0).$

We denote the straight lines
$$L_t=\{u_a(t)+ \lambda E_a(t), \lambda \in \mathbb R \} , \ L_s=\{u_b(s)+\mu E_b(s), \mu \in \mathbb R \}.$$

The tangent spaces at $a$ and $b$ are spanned by the corresponding partial derivatives, which are linearly independent due to regularity:
$$T_a(S)=span \{\partial_t u(t_0,\lambda_0),e(t_0)\},$$
$$T_b(S)=span \{\partial_{t} u(s_0,\mu_0),e(s_0)\}.$$
We know that the function
$$\lambda \to |u(t_0,\lambda)-u(s_0,\mu)|^2$$
attains minimum at $\lambda=\lambda_0, \mu=\mu_0.$ Therefore, the partial derivatives vanish at $(t_0,\lambda_0).$

Differentiation in $\lambda$ at
$t=t_0, \lambda=\lambda_0$ yields
$$\langle e(t_0),u(t_0,\lambda_0)-u(s_0,\mu_0) \rangle =\langle e(t_0),a-b \rangle =0.$$
Analogously, differentiation in $\mu$ gives
\begin{equation}\label{E:1}
\langle e(s_0),a-b \rangle =0.
\end{equation}

For any pair $L_t,L_s$ of the constituting the surface $S$ straight lines, denote
$a(t,s),b(t,s)$ the points
$$a(t,s)=u(t)+\lambda(t,s)e(t), \ b(t,s)+\mu(t,s) e(s),$$
belonging to the lines $L_t,L_s$ correspondingly, at which  the distance between the
lines is attained:
$$d(t,s)=dist(L_t,L_s)=|a(t,s)-b(t,s)|.$$
The coefficients $\lambda(t,s), \mu(t,s)$ can be found from the orthogonality conditions
$$ \langle a(t,s)-b(t,s), e(t) \rangle =0, \ \langle a(t,s)-b(t,s), e(s) \rangle =0.$$
The solutions of the corresponding linear system are
$$\lambda(t,s)=\frac{-\langle e(t),e(s) \rangle \langle u(t)-u(s), e(s) \rangle +\langle u(t)-u(s),e(s) \rangle }{1- \langle e(t),e(s) \rangle^2},$$
$$\mu(t,s)=\frac{\langle e(t),e(s) \rangle \langle u(t)-u(s), e(s)\rangle - \langle u(t)-u(s),e(t) \rangle}{1- \langle e(t),e(s) \rangle ^2}.$$
The denominator is different from zero as the lines $L_t, L_s$ are not parallel by the
condition and hence $1- \langle e(t),  e(s) \rangle \neq 0.$

The above formulas show that the functions $\lambda (t,s), \mu (t,s)$ are differentiable
at  the point $(t_0,s_0).$

Since the distance function $d(t,s)$ attains its maximum at $t_0,s_0$ we have
$$\partial_t d(t_0,s_0)=<a^{\prime}_t (t_0,s_0)-b^{\prime}_t(t_0,s_0),a-b>=0,$$
$$\partial_s d(t_0,s_0)=<a^{\prime}_s (t_0,s_0)-b^{\prime}_s(t_0,s_0),a-b>=0,$$
or
$$ \langle u^{\prime}(t_0)+\lambda^{\prime}(t_0,s_0)e(t_0)+\lambda_0 e^{\prime}(t_0), a-b \rangle =0,$$
$$\langle u^{\prime}(s_0)+\mu^{\prime}(t_0,s_0)e(s_0)+\mu_0 e^{\prime}(s_0), a-b \rangle =0.$$
Since $a-b$ is orthogonal to $e(t_0)$ and $e(s_0),$ we obtain:
$$\langle u^{\prime}(t_0)+\lambda_0 e^{\prime}(t_0), a-b \rangle =0,$$
$$\langle u^{\prime}(s_0)+\mu_0 e^{\prime}(s_0), a-b \rangle =0.$$
Therefore the vector $a-b$ is orthogonal to the vectors
$(\partial_{\lambda}u)(t_0,\lambda_0)$ and to $(\partial_{t})u(t_0,\lambda_0)$ which span
the tangent plane $T_a(S).$ Thus,
$$a-b \perp T_a (S).$$
Analogously,
$$a-b \perp T_b(S).$$
That means that the points $a$ and $b$ are antipodal.

\subsection{End of the proof of Theorem \ref{T:Main1}}

\noindent
{\bf The "if" part.}

Notice, that the "if" statement holds in any dimension $d.$
Suppose that $S$ is a harmonic cone with a vertex $a.$ This means that  there exists a
nonzero harmonic homogeneous polynomial (solid harmonic) $h$ such that
$$h(a+x)=0, \ \forall x \in S.$$
By shifting, we  can assume $a=0.$

Define
$$\varphi_{\lambda}(x) =\int_{|\omega|=1} e^{i \lambda  <x,\omega>}h(\omega)dA(\omega).$$
Then
$$\Delta \varphi_{\lambda}=-\lambda^2\varphi_{\lambda}.$$

Now fix $x_0 \in \mathbb R^d \setminus 0$ such that $h(x_0)=0.$  Denote $SO_x(d)$ the group
of orthogonal transformations $\rho \in SO(d)$ of $\mathbb R^d$ such that $\rho(x_0)=x_0.$
Then
$$\varphi_{\lambda}(x_0)=\varphi_{\lambda}(\rho (x_0))=\int_{|\omega|=1} e^{i \lambda
<\rho ( x_0) ,\omega>}h(\omega)dA(\omega)=\int_{|\omega|=1} e^{i \lambda <x_0, \rho^{-1}(
\omega)>}h(\omega)dA(\omega).$$ Change of variables $\omega^{\prime}=\rho^{-1}(\omega)$ leads
to
$$
\varphi_{\lambda}(x_0)= \int_{|\omega^{\prime}|=1} e^{i \lambda
<x,\omega^{\prime}>}h(\rho \omega^{\prime})dA(\omega^{\prime}).$$
Integrating the equality in
$\omega^{\prime}$ against normalized Haar measure $d\rho$ on $SO(d)$ yields
\begin{equation}\label{E:varphi(x)}
\varphi_{\lambda}(x_0)=\int_{|\omega^{\prime}|=1}e^{i \lambda <x_0,\omega^{\prime}>} \tilde
h(\omega^{\prime})d\omega^{\prime},
\end{equation}
where $\tilde h (\omega^{prime})$ is the average
$$\tilde h(\omega^{\prime})=\int\limits_{\rho \in SO_x(d)}h(\rho \ \omega^{\prime})d\rho.$$
The function $\tilde h(\omega^{\prime})$ is a spherical harmonic, invariant under
rotations $\rho \in SO (d),$ preserving $x_0$, and therefore it is proportional to the
zonal harmonic $Z_{x_0}$ (\cite{SW}) with the pole $x_0/|x_0|,$ of the same degree as $\tilde h:$
\begin{equation}\label{E:zonal}
\tilde h= c Z_{x_0}.
\end{equation}
However,
$$\frac{1}{|x|^{deg h}} \tilde h(x)= \tilde h(\frac{x}{|x|})=h(\frac{x}{|x|})=0,$$
because $\rho \ x=x$ and $h(x)=0$ and $h$ is homogeneous.
On the other hand the value of
the zonal harmonic at its pole is
$$Z_x(\frac{x}{|x|})= \alpha \Omega_{d-1}^{-1},$$
where $\alpha $ is the dimension of the space of spherical harmonics of degree $deg h$
and $\Omega_{d-1}$ is the area of the unit sphere in $\mathbb R^d.$ (\cite{SW}, Corollary
2.9), Therefore, we have form (\ref{E:zonal}):
$$c \alpha \Omega_{d-1}^{-1}=0$$
and $c=0.$  Then (\ref{E:zonal} implies  $\tilde h \equiv 0$ and then $\varphi_{\lambda}(x_0)=0$ because of (\ref{E:varphi(x)}.

Thus, we have proven $\varphi_{\lambda}(x_0)=0$ whenever $h(x_0)=0$ and hence the harmonic cone $h^{-1}(0)$ is a common nodal set for a nontrivial Paley-Wiener family of eigenfunctions.

\noindent
{\bf The "only if"  part}

We assume that an irreducible real analytically ruled hypersurface $S \subset \mathbb R^3,$  without parallel generating lines, is contained
in the common zero set of a  Paley-Wiener family of eigenfunctions.
We need to prove that $S$ is a cone.

We start with the case when the  foliation $\{L_t\}$  of $S$ is everywhere regular. In
particular, it is regular at the extremal points $a,b$ at which the distance function
$d(t,s)$ attains its maximum. Then the points $a$ and $b$ are antipodal by Lemma
\ref{L:anti}, and then Theorem \ref{T:antipodal} implies that $S$ is an injectivity set.
By Proposition \ref{P:observations}, this contradicts to the assumption that $S$ is the
common nodal set for Paley-Wiener family of eigenfunctions.

Therefore $S$ has at least one singular point, say, $a.$ By Corollary \ref{C:harm} of
Theorem \ref{T:conical}  $a$ is a conical point. This means that $a$ belongs to an open
family of lines $\{L_t\}.$ Since $S$ is irreducible, the base curve $\gamma$ that
parametrizes the family $L_t$ is real analytic and connected. Therefore, all lines $L_t$
pass through $a$ and therefore $S$ is a cone with the vertex $a.$

Moreover, $S$ is a harmonic cone. Indeed, we know from Proposition \ref{P:observations} that
there exists a nonzero harmonic polynomial $H$ such that $S \subset H^{-1}(0).$ Since $S$
is a cone with the vertex $a$ we have
$$H(a+\lambda (x-a))=0$$
for all $x \in S$ and $\lambda \in \mathbb R.$ Therefore, if
$H(a+u)=\sum_{j=0}^{N}H_j(u)$ is the homogeneous decomposition, then $H_j(x-a)=0,
j=0,...,N$ and it remains to note that all $H_j$ are harmonic and homogeneous. Then $a+S
\subset h^{-1}(0),$ where $h$ can be taken  any nonzero polynomial $H_j.$   Theorem
\ref{T:Main1} is proved.

\section{Reducible case. Proof of Theorem \ref {T:mainmain1}}

 Now we turn to  the proof of more general Theorem \ref{T:mainmain1} where we do not
 assume that the base curve $\gamma$ of the ruled surface $S$ is connected.

In general situation, $S$ decomposes into irreducible components:
$$S=\cup_{j=1}^M S_j,$$
where each $S_j$ is a real analytically ruled surface with a real analytic closed
connected base curve $\gamma_j.$ So, the ruled surface $S$ is parametrized by the base
curve $$\gamma=\gamma_1 \cup ...\cup \gamma_M.$$

Each surface $S_j$ satisfies all the conditions of Theorem \ref{T:Main1} and therefore is
a harmonic cone with a vertex $a_j \in S_j.$ All we need now is to prove the additional
properties of the decomposition of $S$ into union of cones , claimed in Theorem
\ref{T:mainmain1}.

We will start with proving that the cones pairwise meet.
\begin{Lemma}\label{L:disjoint} If there are $i, j$
such that $S_i \cap S_j =\emptyset$ then $S$ is an injectivity set.
\end{Lemma}

\noindent
{\bf Proof}  Assume that $S$ fails to be an injectivity set. Since $S_i$ and $S_j$ do not
meet, any two generating lines $L_a, \ a \in \gamma_i$ and $L_b ,\ b \in S\gamma_j,$ are
disjoint and $dist (L_a,L_b) >0.$

Since there are no parallel generating lines, the function $(a,b) \to dist(L_a, L_b)$ is
continuous and attains its {\it minimum}.  Let $a_0 \in S_i, b_0 \in S_j$ are the points
where the minimal distance between the generating lines is realized:
$$|a_0-b_0|=\min_{a \in S_i, b \in S_j }  dist (L_a, L_b)>0 .$$
The two cases are possible:
\begin{enumerate}
\item $a_0$ and $b_0$ are regular points of the foliation $S=\cup_{a \in S} L_a.$
\item One of the points $a_0, b_0$ is a singular point.
\end{enumerate}

Let $a_0=u(t_0,\lambda_0), \ b_0=u(s_0,\mu_0).$

In the case 1, the equations:
$$\frac{\partial |u(t,\lambda)-u(s,\mu)|} {\partial t}(t_0,\lambda_0,s_0,\mu_0)=
\frac{\partial |u(t,\lambda)-u(s,\mu)|} {\partial s}(t_0,\lambda_0, s_0, \mu_0)=0,$$
$$\frac{\partial |u(t,\lambda)-u(s,\mu)|} {\partial \lambda}(t_0,\lambda_0,s_0,\mu_0)=
\frac{\partial |u(t,\lambda)-u(s,\mu)|} {\partial s}(t_0,\lambda_0, s_0,\mu_0)=0$$ yield
that the vector $a_0-b_0$ is orthogonal to the tangent spaces $T_{a_0}(S), \ T_{b_0}(S).$
In other words, $a_0$ and $b_0$ are antipodal points. By Lemma \ref{L:antipodal}, $S$ is
an injectivity set. Contradiction.

Consider now the case 2, i.e., assume that one of the extremal points, say, $a_0$ is
singular. Since  $S$ is not an injectivity set, $a_0$ is a conical point, due to Theorem
\ref{T:conical}. The ruled surface $S_i$ has the real analytic connected base curve
$\gamma_i$ hence $S_i$ is a cone with the vertex $a_0.$

Now, the straight lines $L_{t_0}\subset S_i$ and $L_{s_0} \subset S_j$ are the
closest generating lines belonging to $S_i$ and $S_j$ correspondingly. Since $a_0 \in
L_{t_0}, b_0 \in L_{s_0}$ are the closest points, we have
$$L_{t_0}, L_{s_0} \perp [a_0,b_0].  $$

However, since $S_i$ is the cone with the vertex $a_0,$ all the straight lines $L_t$
generating $S_i$ all pass through $a_0.$ If $L_{t}$ is not orthogonal to $[a_0,b_0]$ then
$$dist (L_t,L_{s_0}) < |a_0-b_0|=dist(L_{t_0},L_{s_0})$$ which is impossible.

Therefore, for all generating lines $L_t \subset S_i$ we have
$$L_t \perp [a_0,b_0]$$ and
hence $L_t \subset \Pi,$ where $\Pi$ is the plane passing through $a_0$ and orthogonal to
$[a_0,b_0].$ Then $S_i$ coincides with the plane $\Pi$ and $S_i=\Pi$ can be viewed as a
line foliation, regular at $a_0.$ If the second extremal point $b_0$ is regular for the
given foliation $\{L_t\}$ then both points $a_0,b_0$ are regular antipodal points and $S$
is an injectivity set. If $b_0$ is a conical point, then the same argument with closest
generating lines shows that $S_j$ is a plane. Then again $a_0,b_0$ are regular antipodal
points and $S$ is an injectivity sets. Lemma is proved.

Now we will prove that the cones intersect transversally.

\begin{Lemma} \label{L:tangent} If some $S_i$ and $S_j$ are tangent at a point
$a$ which is not a vertex of any cone $S_i,S_j$ then $S$ is an injectivity (not nodal) set.
\end{Lemma}

\noindent
{\bf Proof} We saw in the proof Theorem \ref{T:conical} that if $a$ is not a vertex of the cone $S_i$ then it is either the point of real analyticity
or a  point of differentiability , which is a singular point of the line foliation and corresponding to the case of odd $m$ in the parametrization (\ref{E:canonical}).
The same is true for the cone $S_j.$

After a suitable translation and rotation, we can make $a=0$ and
$$T_a(S_i)=T_a(S_j)=\{x_3=0\}.$$

The representation (\ref{E:Puis}) shows that the surfaces $S_i,S_j$ are defined near $a=0$ as the graphs:

\begin{eqnarray*}
S_i:x_3 &=& z_i(x_1,x_2), \\
S_j:x_3 &=& z_j(x_2,x_2),
\end{eqnarray*}

where
$$z_i(x_1,x_2)=o(r), \ z_j(x_1,x_2)=o(r), \ r=\sqrt{x_1^2+x_2^2} \to 0.$$
Moreover, by the construction, these functions are algebraic and for some odd integers $m,n$ the functions
$$ z_i(x_1^m,x_2), z_j(x_1^n,x_2)$$
are real analytic.

If $S$ is not an injectivity set, then  due to Proposition \ref{P:observations}, there exists the nonzero harmonic polynomial $H$ vanishing on $S$
(Proposition \ref{P:observations}).
Since $H=0$ on $S_i=\{x_3-z_i(x_1,x_2) =0,\}$
the polynomial
$$H(x_1^{mn},x_2,x_3)=0 \ \mbox{ whenever}  \ \rho_i(x):=x_3-z_i(x_1^{mn},x_2)=0.$$
The function $\rho$ is real analytic and $\nabla \rho \neq 0$ hence the polynomial $H$ is divisible by $\rho$ which means that
$$H(x_1^{mn},x_2,x_3)=(x_3- z_i(x_1^{mn},x_2)) R(x_1,x_2,x_3),$$
where $R$ is real analytic near $0.$

Since $S_i$ and $S_j$ can coincide only on a nowhere dense subset, and $H=0$ on $S_j,$ the function $R$ must vanish on the surface
$\rho_j(x):=x_3-z_j(x_1^{mn},x_2)=0.$
Further, since both functions $H$ and $\rho_j$ are real analytic and $\nabla \rho_j \neq 0,$ the function $R$ is divisible by $\rho_j,$ meaning that
$$R=\rho_j G,$$
where the function $G$ is  real analytic near $0.$

Finally, returning to $x_1$ instead of $x_1^{mn}$   we have
$$H(x)=(x_3-z_i(x_2,x_3))(x_3-z_j(x_1,x_2))G(x^{\frac{1}{mn}},x_2,x_3).$$

Decompose
$$G(x^{\frac{1}{mn}},x_2,x_3)=\sum_{\alpha, \beta,\gamma \geq 0}x_1^{\frac{\alpha}{mn}}x_2^{\beta} x_3^{\gamma}$$
and let $G_0$ be the sum of the terms with the minimal homogeneity degree $$\frac{\alpha_0}{mn}+\beta_0+\gamma_0.$$

If
$$H=H_k+H_{k+1}+...+H_N, H_k \neq 0,$$
is the homogeneous decomposition for $H,$ then since $z_i, z_j=o(r), r \to 0$ we have for the minimal degree homogeneous term:
$$H_k(x)=x_3^2 G_0(x).$$
Thus,
$$H(x_1,x_2,0)=0.$$
Notice that $G_0$ is a polynomial with respect to $x_1^{\frac{1}{mn}}, x_2, x_3.$
Therefore, differentiation in $x_3$ yiedls:
$\partial_{x_3} H(x)= 2x_3 G_0(x)+ x_3^2 \partial_{x_3}G_0(x)$ and hence
$$\partial_{x_3}H(x_1, x_2,0)=0.$$
However, the polynomial $H_k$ is harmonic and satisfies the overdetermined Dirichlet-Neumann conditions on the plane $x_3=0.$
This implies $H_k \equiv 0.$ This contradiction completes the proof.

\subsection{End of the proof of Theorem \ref{T:mainmain1}}
First of all, according to Theorem \ref{T:Main1}, each irreducible ruled component of $S$
is a harmonic cone and therefore, $S$ is the union of harmonic cones, $S=\cup_{j=1}^N
S_j.$

Moreover, the vertices are the only singular points of the cones $S_i.$ The cones
$S_i$ are real analytic everywhere except, maybe, for the vertex. If $S_i$ is
differentiable at the vertex then $S_i$ is a plane and, of course, is real analytic
everywhere.

Further, Lemma \ref{L:disjoint} implies that $S_i \cap S_j \neq \emptyset$ for any $i
\neq j,$ since otherwise $S$ is an injectivity set. In turn, Lemma \ref{L:tangent} says that $S_i \neq S_j$ is transversal. The intersection
$S_i \cap S_j$ is either 0-dimensional (discrete) or one-dimensional. In the latter case
the intersection is a curve.

In the case when $S_i \cap S_j$ is discrete, then since $S_i, S_j$ are two-dimensional,
any  point $a \in S_i \cap S_j,$ at which $S_i$ and $S_j$ are differentiable, must be a
tangency point, which is not the case. Therefore, $a$ must be singular for either cone
$S_i,S_j$  and hence is a  vertex of one of them.  Theorem \ref{T:mainmain1} is proved.

\section{ Coxeter systems of planes. Proof of Theorem \ref{T:mainmain2} }

Theorem \ref{T:Main1} asserts that $S$ is a cone. The only cone which has no
differentiable singularities is a plane. Therefore, if $S$ in Theorem \ref{T:Main1} is differentiable surface then $S$ is a plane.

Then Theorem{T:mainmain2} follows from
\begin{Lemma}\label{L:union} Any finite union $S$ of hyperplanes in $\mathbb R^d$ is an injectivity
set unless $S$
can be completed to a Coxeter system.
\end{Lemma}

\noindent
{\bf Proof}

We will give the proof for the case  $d=3$ which is under consideration in this article.

Let $$S=\cup_{i=1}^{N} \Pi_i$$
where $\Pi_i$ are the hyperplanes.  Suppose that $S$ fails to
be an injectivity set. Then there exists a nonzero   function  $f \in C_{comp}(\mathbb
R^3)$ such that $Rf(x, t)=0, \ t>0,$ for all $x \in S.$  It is known \cite{CH}, v.II,
that then $f$ is odd with respect to reflections around each plane $\Pi_i.$

 Denote $W_{\Pi_1,...,\Pi_N}$ the group generated by the reflections around the planes $\Pi_1,...,\Pi_N.$

Now we are going to use the additional information about existence of nonzero harmonic
polynomial vanishing on $S$ (Proposition \ref{P:observations}), which rules out, due to
Maximal Modulus Principle, the possibility for the action of the group
$W_{\Pi_1,...,\Pi_N}$ to have compact fundamental domain.

If $N=2$ then the angle between $\Pi_1$ and $\Pi_2$ must be a rational multiple of
 $\pi$ since otherwise $$\cup_{w \in W_{\Pi_1,\Pi_2}} w (\Pi_1)  \cup \cup_{w \in W_{\Pi_1,\Pi_2}} w (\Pi_2)$$ is dense in $\mathbb R^3$ and then $f=0$
 identically because $f$ vanishes on each $\Pi_1, \Pi_2.$ Therefore $S$ is a subsystem of the Coxeter system
generated by the planes $\Pi_1,\Pi_2.$

Let $N \geq 3.$  The following cases are possible:
\begin{enumerate}
\item all the planes $\Pi_i, i=1,...,N,$ have a common point,
\item there are two parallel planes $\Pi_{i_1},\Pi_{i_2}$,
\item there are  three planes $\Pi_{i_1},\Pi_{i_2},\Pi_{i_3}$  that bound  a right  triangular prism,
\item $N \geq 4$ and there are four planes $\Pi_{i_1},\Pi_{i_2},\Pi_{i_3}, \Pi_{i_4}$ that bound a bounded
simplex.
\end{enumerate}

In the first case, the reflection group $W$ generated by the planes $\Pi_i$ must be finite,
since otherwise $ \cup_{w \in W_{\Pi_1,...,\Pi_N}} w(S)$ is dense in $\mathbb R^3$ and then $f=0.$
Therefore, in the first case $S$ can be included in a Coxeter system of planes.

The second case is impossible, since $supp f,$ being symmetric
both with respect to $\Pi_{i_1}$ and $\Pi_{i_2},$  must be unbounded, which is not the case.

In the third case, the normal vectors $\nu_1,\nu_2,\nu_3$ of the corresponding planes are
linearly dependent and span a plane $P$ orthogonal to all $P_{i_j}, j=1,2,3.$ For any $b
\in \mathbb R^3$ the intersection $(P+b) \cap (P_{i_1} \cup \Pi_{i_2} \cup \Pi{i_3})$  is
three lines $L_1,L_2,L_3$ in the 2-plane $P+b,$ bounded a triangle.

The restriction $f\vert_{P+b}$ can be regarded as a compactly supported function defined in
$\mathbb R^2,$  and this function is odd-symmetric with respect to the lines
$L_1,L_2,L_3.$

In particular, it has zero spherical means on the lines. As it was proven in Proposition
\ref{P:observations} , if $f$ is not identically zero on $P+b$ then there is a nonzero
harmonic polynomial vanishing on $L_1 \cup L_2 \cup L_3$ which is impossible due to
Maximum Modulus Principle since the union contains a bounded contour. Therefore, $f=0$ on
$P+b$ and then $f=0$ everywhere as $b$ is arbitrary. Thus, the third case is ruled out as
well.

Also, the fourth case is impossible, since if $f$ is not zero then  we again have
contradiction with existence of a nonzero harmonic polynomial vanishing on $S,$   as in
the previous case. Lemma is proved.

\noindent
{\bf Proof of Theorem \ref{T:mainmain2}}  Since any two-dimensional cone in $\mathbb R^3,$ which is a differentiable surface,
is a two-dimensional plane,  Theorem \ref{T:mainmain1} implies that  the surface
$S$ in Theorem \ref{T:mainmain2} is a finite union of 2-planes and hence is
a Coxeter system of planes, due to Lemma \ref{L:union}.

\section{Proof of Theorem  \ref{T:main_conv_nodal} (the case of convexly supported generating function)}

\subsubsection{Lemmas}

We are given a nonzero function $f \in C_{comp}(\mathbb R^3).$ Consider the set
$$N_f=\{x \in \mathbb R^3: Rf (x,t)=0, \ \forall t>0 \}.$$  By Proposition
\ref{P:observations}, the set $N_f$ represents as
$$N_f=S \cup V,$$
where $S$ is either empty or an algebraic hypersurface
$$S=Q^{-1}(0),$$
where $Q$ is
a polynomial, dividing a nonzero harmonic polynomial $H.$ We assume $S \neq \emptyset.$

Denote $\Gamma$ the outer boundary of $supp f.$ By the condition, $\Gamma$ is strictly
convex real analytic closed hypersurface. Theorem \ref{T:ruled} yields  that the
observation surface $S$ is foliated into straight lines, each of which intersects
orthogonally, at two points, the strictly convex surface $\Gamma.$

The surfaces $\Gamma$ and $S$ intersect orthogonally. The intersection
$$\gamma: \Gamma \cap S$$
is a curve, smooth at all points $a \in \gamma$ at which $S$ is smooth.

\begin{Lemma} \label{L:Sj} The surface $S$ is a real analytically ruled surface.
\end{Lemma}

\noindent
{\bf Proof}
Denote
$$\gamma =\Gamma \cap S.$$
Pick a point $ a \in \gamma.$  Let $T_a(\Gamma)$ be the tangent plane. Applying
translation and rotation, one can assume that $a=0$ and $$T_a (\Gamma)=\{x_3=0\}.$$

The projection
$$\pi : T_a(\Gamma) \mapsto \Gamma$$
along the normals to $\Gamma$ is well defined in a neighborhood $$U\subset T_a(\Gamma)$$
of $a.$

Since $\Gamma$ is real analytic, the normal field to $\Gamma$ is real analytic as well and hence
$\pi$ is real analytic diffeomorphism near $a=0.$
Also, $\pi(U \cap S)$ is an open neighborhood of $a$ in  $\gamma.$

It is easy to understand that the polynomial $Q$ is not identically zero on $T_a(S)$
since $S=Q^{-1}(0)$ is transversal to $T_a(S)$ near $a.$ Therefore, the intersection
$$C:= T_a(\Gamma) \cap S $$
is an open algebraic curve in the plane $T_a(\Gamma)=\{z=0\},$ defined by the equation
$C=\{Q(x,y,0)=0\}.$

Then we  use  Puiseux theorem (\cite{L}, Ch.II, 9.6; \cite{W}, Thm. 2.1.1; \cite{C},
Ch.2,p. 3-11) which claims that each branch $C_i$ of $C$ is parametrized either by
$$I \ni t \mapsto (0,t,0),$$
or by
$$I \ni t \mapsto (t^m, \alpha_i(t),0)$$
where $I$ is an open interval (which can be taken $I=(-1,1)$), $m$ is natural and $\alpha_i(t)$ is a
real analytic function.

Then   $\gamma$ decomposes, near $a,$ into the union of the curves $\gamma_i=\pi(C_i)$
and each $\gamma_i$ is the image $\gamma_i=u(I)$ where the mapping
$$ I \ni  \mapsto u_i(t)= \pi(t^m,B_i(t),0)$$ is real analytic, because $\pi$ is so.
By Corollary \ref{C:harm} of Theorem \ref{T:conical}, the ruled surface
$$S_i=\{u_i(t)+\lambda \nu(u_i(t)), t \in I, \lambda \in \mathbb R\},$$ where $\nu$ is
unit normal vector to $\Gamma,$ is real analytically ruled surface.

\begin{Lemma} \label{L:midpoint} Let $a$ be the vertex of the cone $C_i.$ Let $\gamma_i$
be a connected closed subarc of $C_i \cap \Gamma$ where $\Gamma$ is the outer boundary of
$supp f.$  Then the distance $|x-a|$ from $a$ to an arbitrary  point $x \in \gamma_i$ is
constant.
\end{Lemma}

\noindent
{\bf Proof}  Consider the parametrization $u(t,\lambda)=u(t)+\lambda e(t), t \in I, $ of
the cone $C_i.$ The mapping $t \mapsto u(t)$ parametrizes the curve $\gamma_i =C_i \cap
\Gamma.$ Consider the distance function
$$d(t)=|a-u(t)|^2. $$
Then
$$d^{\prime}(t)=(a-u(t),u^{\prime}(t)).$$
Since $a$ is the vertex of $C_i,$ it belongs to any line $L_t.$ Therefore
$a=u(t)+\lambda(t) e(t)$ and hence
$$d^{\prime}(t)=(a-u(t), u^{\prime}(t))=\lambda(t)(e(t),u^{\prime}(t))=0,$$
because $u^{\prime}(t)$ is tangent to $\Gamma,$ $e(t)$ is the directional vector of the
Line $L_t$ and $L_t$ is orthogonal to $\Gamma,$ as stated in Theorem \ref{T:ruled}.

\begin{Lemma} \label{L:outside} If two cones $C_i, \ C_j$ meet outside of $supp f$
then they have a common vertex and hence the union $C_i \cup C_j$ is itself a cone.
\end{Lemma}

\noindent
{\bf Proof} The cones $C_i, C_j$ consist of straight lines orthogonal to the outer
boundary $\Gamma$ of $supp f.$ Also, $\Gamma$ is a real analytic strictly convex surface.
 If $C_i$ meet $C_j$ in the exterior of
$\Gamma$ then $C_i$ and $C_j$ share a ruling straight line $L$ passing through a common
point of the two cones and orthogonal to $\Gamma.$ The vertices of both cones $C_i$ and
$C_j$ belong to $L.$ The common line $L$ meets the convex surface $\Gamma$ at two points
$b^{+}, b^{-}:$
$$\{ b^+,b^{-} \}= L \cap \Gamma.$$

Let $\gamma_i$ and $\gamma_j$ be the connected closed subarcs of  the smooth curves $C_i
\cap \Gamma$ and $C_j \cap \Gamma,$ correspondingly, containing the point $b^+.$

Then $\gamma_i, \gamma_j$ are smooth closed curves on $\Gamma,$ sharing the common point
$b^{+} \in \gamma_i \cap \gamma_j.$

Suppose that $\gamma_i$ and $\gamma_j$ are tangent at $b^+$ and let $\tau$ be the common
tangent vector at $b^+.$ Since the tangent planes of the cones $C_i$ and $C_j$ coincide:
$$T_{b^+} (C_i)=T_{b^+} (C_j)= span \{L,e\},$$
the two cones are tangent. However, this is impossible due to Lemma \ref{L:tangent}.

Thus, the two closed curves $\gamma_i$ and $\gamma_j$ intersect at $b^+$ transversally.
Then they must intersect in at least one more point, $c \in \Gamma.$  Then both cones
$C_i, C-j$ contain the straight line $L_c$ intersecting $\Gamma$ orthogonally at the
point $c.$ The two cases are possible:
\begin{enumerate}
\item $c \neq b^{-}.$
\item $c=b^{-}.$
\end{enumerate}

In the case 1 , the straight lines $L$ and $L_c$ are different. Both of them belong to
the cones $C_i$ and $C_j$ and hence the intersection of the two lines $L \cap L_c$ is
just a single point which is the vertex of both $C_i$ and $C_j.$ Thus, $C_i$ and $C_j$
share the vertex and Lemma is proved in this case.

In the case 2 the two straight lines coincide, $L=L_c,$ as they both pass through the
points $b^+$ and $b^-=c.$ Let  $a_i,a_j$ be the vertices of the the cones $C_i,C_j$
correspondingly. By Lemma \ref{L:midpoint}, the distance $|x-a_i|$ is constant on
$\gamma_i$. Since $b^+, b^{-} \in \gamma_i$, we have
$$|b^+ - a_i|=|b^- - a_i|.$$
The three points $a, b^+, b^-$ belong to the same line $L$ and therefore, $a$ is the
midpoint:
$$a_i=\frac{1}{2}(b^+ + b^-).$$
The same can be repeated for $\gamma_j$ and then we obtain
$$a_j=\frac{1}{2}(b^+  + b^{-}).$$
Thus, $a_i=a_j$ and the statement of Lemma is true in the case 2 as well.

\begin{Lemma}\label{L:0-dim}  Suppose that $S_i \cap S_j$ is 0-dimensional. Then
\begin{enumerate}
\item
$S_i \cap S_j \subset \{c_i,c_j\},$ where $c_i,c_j$ are the vertices of the cones
$S_i,S_j$ correspondingly.
\item If $S_i \cap S_j =\{c_i,c_j\}$ then $c_i=c_j.$
\end{enumerate}
\end{Lemma}

\noindent
{\bf Proof} We know that $S_i$ and $S_j$ are differentiable everywhere except maybe at
the vertices. If $a \in S_i \cap S_j$ and $a \neq c_i, a \neq c_j,$ then $a$ is the point
of smoothness for both $S_i$ and $S_j$  and hence the cones $S_i, S_j$ cannot intersect
at $a$ transversally since in this case the intersection $S_i \cap S_j$ must be
one-dimensional. Therefore, $S_i$ and $S_j$ are tangent at $a.$ This possibility is ruled
out by Lemma \ref{L:tangent}. This proves the statement 1.

If $S_i \cap S_j=\{c_i,c_j\}$ and $c_i \neq c_j$ then both cones $S_i$ and $S_j$ contain
the straight line passing through the vertices $c_i$ and $c_j.$ This contradicts to the
assumption that the intersection is 0-dimensional.

\begin{Lemma}\label{L: 1-dim}  If $S_i \cap S_j$ is one-dimensional then the cones $S_i$ and $S_j$ share
the vertex so that $S_i \cup S_j$ is a cone.
\end{Lemma}

\noindent
{\bf Proof} Let $\gamma=S_i \cap S_j.$ If the curve $\gamma$ is unbounded, then $S_i$ and
$S_j$ intersect outside of $\Gamma$ and by Lemma \ref{L:outside} $S_i$ and $S_j$ have a
common vertex. Otherwise, $\gamma$ is a bounded curve. It is also closed as it is
algebraic. Then $\gamma$ bounds two -dimensional domains $D_i$ and $D_j$ on the surfaces
$S_i,S_j$ correspondingly. Therefore, $S_i \cap S_j$ contain a cycle $D_i \cup D_j.$
However, it is impossible due to Maximum Modulus Principle, since there exists a nonzero
harmonic polynomial $H$ vanishing on $S_i \cup S_j.$

\begin{Corollary} \label{C:corr} If $S_i$ and $S_j$ have different vertices,
$c_i \neq c_j,$ then $S_i \cap S_j$ consists of a single point, which is
either $c_i$ or $c_j$.
\end{Corollary}

\noindent
{\bf Proof} The intersection  $S_i \cap S_j$ is discrete (0-dimensional) since otherwise the cones $S_i, S_j$ have equal vertices, by Lemma \ref{L: 1-dim}.
Then Lemma \ref{L:0-dim} says the intersection coincides with one of the vertices.

\subsubsection{ End of the proof of Theorem \ref{T:main_conv_nodal}}

Let us group all the cones $S_i$ whose vertices coincide. The union of  such cones is
again a cone  and hence the union $S$ can be regrouped in the union
$$S=C_1 \cup...\cup C_P$$
of cones $C_i$ with pairwise different vertices $b_i.$. Each $C_i$ is the union of the
cones $S_j$ with equal vertices. Due to Lemma \ref{L:1-dim}, the pairwise intersections
$C_i \cap C_j, \ i\neq j,$ are 0-dimensional.

First of all , all the cones $C_j$ are harmonic.  Indeed, we know that there is a nonzero
harmonic polynomial $H$ vanishing on $S.$ By translation, we can assume the the vertex
$b_i$ of the cone $C_i$ is $b_i=0.$ Since $C_i$ is a cone, we have
$$H(\lambda x)=0$$
for all $x \in C_i$ and all $\lambda \in \mathbb R.$ If $H=H_0+...+H_N$ is the
homogeneous decomposition, then $H_0(x)+\lambda H_1(x)+...+\lambda^N H_N(x)=0$ and hence
$H_k(x)=0$ for all $k.$ If $h=H_j$ is any nonzero homogeneous polynomial then $h(x)=0$
for all $x \in C_i$ and hence $C_i$ ia a harmonic cone.

Further, we know that for any  $i \neq j$ the intersection $S_i \cap S_j$ is either $c_i$
or $c_j.$ It follows that for the cones $C_i,$ which are unions of groups of $S_j,$ holds
$C_i \cap C_j \subset \{b_i,b_j\}.$ If $C_i \cap C_j=\{b_i,b_j\}$ then both cones $C_i$
and $C_j$ contain the points $b_j \neq b_j$ and hence, the straight line through these
points, which is not the case.

Thus, $C_i \cap C_j$ is a single point,which is a vertex of $C_1$ or $C_2:$
\begin{equation}\label{E:either}
C_i \cap C_j=\{b_i\} \ \mbox{or} \ \{b_j\}.
\end{equation}

\begin{Lemma} \label{L:P<4}
$P \leq 3.$
\end{Lemma}

\noindent
{\bf Proof} Suppose that $P \geq 4.$ Consider the cones $C_1,C_2,C_3,C_4.$ We have
$$C_1 \cap C_2=\{b_1\} \ \mbox{or} \ \{b_2\}.$$
Without loss of generality, we can assume that
$$C_1 \cap C_2=\{b_1\}.$$
Then
$$C_1 \cap C_3=\{b_3\}.$$
Indeed, if $C_1 \cap C_3=\{b_1\}$ then $b_1 \in C_3, b_1 \in C_2$ and  therefore $$b_1
\subset \{b_2,b_3\},$$
 which is impossible because  $b_1,b_2,b_3$ are  all different. For the same reason,
$$C_1 \cap C_4=\{b_4\}.$$
Now,
$$C_2 \cap C_3=\{b_2\},$$
because otherwise $C_2 \cap C_3 =\{b_3\}$ and then $b_3 \in C_2, b_3 \in C_1$ and
therefore $$b_3 \in \{b_1,b_2\}$$ which is not the case.

Now consider the intersection of $C_2$ and $C_4:$
$$C_2 \cap C_4=\{b_2\} \ \mbox{or} \ \{b_4\}.$$
If $C_2 \cap C_4=\{b_2\}$ then we have
$$b_2 \in C_4, b_2 \in C_3$$ and therefore $$b_2 \in \{b_3,b_4\}$$ which is not the case.
If, alternatively,  $C_2 \cap C_4=\{b_4\},$ then we have $b_4 \in C_2, b_4 \in C_1$ and
therefore $$b_4 \in \{b_1,b_2\},$$ which is not the case.  Thus, neither option is
possible. Thus, $P \leq 3.$ Lemma is proved.

Let us continue the proof of Theorem \ref{T:main_conv_nodal}.

If $P=1$ then $S=C_1$ is a cone and, moreover, a harmonic cone. This is the case 1) in
Theorem \ref{T:mainmain1}.

Suppose $P=2$ so that $S=C_1 \cup C_2.$ Formula (\ref{E:either}) leads to the case 2) of
Theorem \ref{T:mainmain2}.

Finally, suppose that $P=3$  and therefore $$S=C_1 \cup _2 \cup C_3.$$
\begin{Lemma}\label{L:two_vertices}
No two cones of $C_1,C_2,C_3$  can have vertices belonging to the third one.
\end{Lemma}

\noindent
{\bf Proof}
Suppose, for example,  that
$$b_1, b_2 \in C_3.$$

We know that $C_1 \cap C_2$ is either $b_1$ or $b_2.$ In the first case we have $b_1 \in
C_2$ and also $b_1 \in C_3.$ Hence
$$b_1 \in C_2 \cap C_3.$$ This implies that either $b_1=b_2$ or $b_1 =b_3.$ Neither is
possible as all the vertices are different.

In the second case we have $b_2 \in C_1$ and also $b_2 \in C_3.$ Then $b_2 \in C_1 \cap
C_3,$ which is either $b_1$ or $b_3$ and we have the same kind of contradiction. Lemma is
proved.

Now we can finish the proof of Theorem \ref{T:main_conv_nodal} in the case $S=C_1 \cup
C_2 \cup C_3.$

We have $C_1 \cap C_2=$ is either $b_1$ or $b_2.$ If
$$C_1 \cap C_2=\{b_1\},$$
then $C_2 \cap C_3$ can be only $b_2$ since otherwise $
b_1,b_3 \in C_2$ which is ruled out by Lemma \ref{L:two_vertices}. Analogously, $C_3 \cap
C_1$ cannot be equal to $b_1$ since then $b_1 \in C_2 \cap C_3$ and hence $b_1$ is either
$b_2$ or $b_3$ which is not the case.

The case $C_1 \cap C_2 =\{b_2\}$ is treated in a similar way. Thus, finally we conclude
that in the case $P=3$ the configuration of the cones is exactly as it is pointed out in
the case 3 of Theorem \ref{T:main_conv_nodal}.

Theorem is proved.

\section{Concluding remarks}

\begin{itemize}
\item
Proving in full Conjecture \ref{C:Conjecture} for ruled surfaces requires proving that the configurations of
cones in Theorem \ref{T:mainmain1} is itself a cone, i.e., the vertices of all the cones $C_i$  coincide.
\item Proving Conjecture \ref{C:Conjecture} in general case requires proving that  common nodal sets for
Paley-Winer families of Laplace eigenfunctions are ruled surfaces. Then one could apply
Theorems \ref{T:Main1} and \ref{T:mainmain1} to pass from ruled surfaces to cones.
\end{itemize}

\subsection*{Acknowledgments}
Part of this work was done when the author spent his sabbatical semester in Spring 2014 at
the  University of Birmingham. The author thanks Department of Mathematics of that
university and Professors Jonathan Bennett and Neal Bez,  for the invitation
and the  hospitality.
The author thanks Dmitry Kerner, Inna Scherbak and Eugenii Shustin   for fruitful discussions on
singularities of ruled surfaces, Plamen Stefanov for pointing out proper references related to the microlocal analysis background of Theorem \ref{T:local_symm}, and Yehonathan Salman for reading the manuscript and making helpful remarks.

\end{document}